\documentclass[english]{amsart}

\usepackage{enumerate}
\usepackage{amsthm,amssymb}
\usepackage[margin=0.90in]{geometry}
\usepackage{tikz}
\usepackage{pgf}
\usepackage{pgfplots}
\usepackage{pgfplotstable}
\usepackage{siunitx,booktabs,orcidlink}
\usepackage{subcaption}
\usepackage{algorithm}%
\usepackage{algorithmicx}%
\usepackage{algpseudocode}%
\usepackage[square]{natbib}%
\setcitestyle{numbers}
\AtBeginEnvironment{proof}{\vspace{-12pt}}

\newtheorem{theorem}{Theorem}
\newtheorem{proposition}[theorem]{Proposition}%

\newtheorem{lemma}[theorem]{Lemma}%
\newtheorem{corollary}[theorem]{Corollary}

\newtheorem{example}{Example}%
\newtheorem{remark}{Remark}%

\newtheorem{definition}{Definition}

\numberwithin{equation}{section}

\DeclareMathOperator{\inte}{int}

\DeclareMathOperator{\cone}{cone}
\DeclareMathOperator{\Hess}{Hess}
\DeclareMathOperator{\grad}{grad}
\DeclareMathOperator{\argmin}{arg\,min}

\DeclareMathOperator{\Proj}{Proj}
\newcommand{\Hy}{{\mathbb H}^n_{\kappa}}

\newcommand{\R}{\mathbb R}

\newcommand{\J}{\textrm{J}}
\newcommand{\K}{K}
\newcommand{\C}{C}

\newcommand{\V}{\mathcal V}
\newcommand{\lng}{\left\langle}
\newcommand{\rng}{\right\rangle}
\newcommand{\lf}{\left}
\newcommand{\rg}{\right}
\newcommand{\p}{\partial}
\newcommand{\ovl}{\overline}

\newcommand{\Lo}{L}
\newcommand{\f}{\frac}
\newcommand{\tp}{^\top}

\definecolor{TolLightBlue}{HTML}{77AADD}
\definecolor{TolLightCyan}{HTML}{99DDFF}
\definecolor{TolLightMint}{HTML}{44BB99}
\definecolor{TolLightPear}{HTML}{BBCC33}
\definecolor{TolLightOlive}{HTML}{AAAA00}
\definecolor{TolLightYellow}{HTML}{EEDD88}
\definecolor{TolLightOrange}{HTML}{EE8866}
\definecolor{TolLightPink}{HTML}{FFAABB}
\definecolor{TolMutedBlue}{HTML}{332288}
\definecolor{TolMutedGreen}{HTML}{117733}
\definecolor{TolMutedPurple}{HTML}{AA4499}

\hypersetup{%
    colorlinks=true,%
    citecolor=TolMutedBlue,%
    linkcolor=TolMutedGreen,%
    urlcolor=TolMutedPurple,%
}

\begin{document}


\title[On projection mappings and the gradient projection method]{On projection mappings and the gradient projection method \\on hyperbolic space forms}

\author[R. Bergmann]{R. Bergmann \orcidlink{0000-0001-8342-7218}}
\address[R. Bergmann]{Department of Mathematical Sciences, Norwegian University of Science and Technology, NO-7491 Trondheim, Norway}
\email{\detokenize{ronny.bergmannn@ntnu.no}}
\urladdr{https://www.ntnu.edu/employees/ronny.bergmann}

\author[O. P. Ferreira]{ O. P. Ferreira \orcidlink{0000-0002-5758-0320}}
\address[O. P. Ferreira]{Institute of Mathematics and Statistics, Federal University of Goi\'as, Campus II, 74690-900, Brazil}
\email{\detokenize{orizon@ufg.br}}
\urladdr{https://orizon.ime.ufg.br}

\author[S. Z. N\'emeth]{S. Z. N\'emeth \orcidlink{0000-0001-9768-4530}}
\address[S. Z. N\'emeth]{School of Mathematics, University of Birmingham, Edgbaston, B15 2TT,  United Kingdom}
\email{\detokenize{s.nemeth@bham.ac.uk}}
\urladdr{https://web.mat.bham.ac.uk/S.Z.Nemeth/}

\author[J. Zhu]{J. Zhu}
\address[J. Zhu]{School of Mathematics, University of Birmingham, Edgbaston, B15 2TT,  United Kingdom}
\email{\detokenize{jxz755@student.bham.ac.uk}}

\begin{abstract}
	This paper presents several new properties of the intrinsic $\kappa$-projection into $\kappa$-hyperbolically convex sets of $\kappa$-hyperbolic
	space forms, along with closed-form formulas for the intrinsic $\kappa$-projection into specific $\kappa$-hyperbolically convex sets. It also discusses the relationship between the intrinsic $\kappa$-projection, the Euclidean orthogonal projection, and the Lorentz projection. These properties lay the groundwork for analyzing the gradient projection method and hold importance in their own right. Additionally, new properties of the gradient projection method to solve constrained optimization problems in $\kappa$-hyperbolic space forms are established, considering both constant and backtracking step sizes in the analysis. It is shown that every accumulation point of the sequence generated by the method for both step sizes is a stationary point for the given problem. Additionally, an iteration complexity bound is provided that upper bounds the number of iterations needed to achieve a suitable measure of stationarity for both step sizes. Finally, the properties of the constrained Fermat-Weber problem are explored, demonstrating that the sequence generated by the gradient projection method converges to its unique solution. Numerical experiments on solving the Fermat-Weber problem are presented, illustrating the theoretical findings and demonstrating the effectiveness of the proposed methods.
\end{abstract}


\date{\today}


\maketitle

\section{Introduction}
The development of numerical methods to address constrained optimization problems heavily relies on understanding the geometric structure of the constraint set. It is widely recognized that Riemannian manifolds constitute significant classes of sets with rich intrinsic geometric characteristics and a solid conceptual foundation that can be explored for this purpose, see for instance \cite{Absil2008, Boumal2020}. The hyperbolic space is one of the most extensively studied and applied models of Riemannian geometry with negative curvature, serving as a powerful tool for addressing mathematical, physical, and engineering questions that Euclidean geometry cannot adequately explain. One of its most famous applications is in the theory of relativity, where it serves as a model  \cite{Faber1983}.  Despite the main focus of this paper not being on issues related to applications, it is worth noting that in recent years, hyperbolic space has attracted increasing attention in modeling practical problems due to its unique geometric properties, which offer numerous computational advantages over classical Euclidean spaces.  As extensively discussed in \cite{Moreira2024}, hyperbolic geometry has emerged as a focal point, demonstrating remarkable efficacy in representing hierarchical data structures with consequences across various domains. Traditional Euclidean optimization methods may face challenges in effectively capturing and exploiting these structures, resulting in suboptimal solutions. In contrast, hyperbolic optimization methods leverage the intrinsic geometry of hyperbolic spaces to navigate and optimize over these structures more efficiently. This capability facilitates more precise data modeling and analysis across a wide range of applications, encompassing machine learning, data mining, network analysis, and beyond.  For instance, in the papers  \cite{chami2021horopca, fan2022nested}, methods for hyperbolic dimensionality reduction are introduced. This constitutes a foundational challenge in machine learning, with applications not only in computer vision but also across diverse domains in engineering the sciences in general. These methods extend Euclidean principal component analysis (PCA) to hyperbolic spaces and better preserve information in the original data compared to previous PCA generalizations, thereby improving distance preservation and downstream classification performance. In \cite{jawanpuria2019low} is proposed an  efficient algorithms for learning low-rank factorizations of hyperbolic embeddings, casting them as manifold  optimization problems. In paper \cite{keller2021hyperbolic}, the use of data from European banking stress tests demonstrates alignment between the latent geometry of financial networks and hyperbolic geometry. This connection facilitates monitoring structural changes in financial networks and distinguishing between systemic relevance and peripheral structural changes through hyperbolic embedding methods. The paper \cite{peng2021hyperbolic} reviews neural components for building hyperbolic deep neural networks and extends prominent deep learning methodologies to hyperbolic space. It examines applications across various machine learning tasks on publicly available datasets, providing insights and delineating future research directions. Additional references addressing this subject include but are not limited to \cite{nickel2018learning, muscoloni2017machine, sharpee2019argument,  tabaghi2020hyperbolic, Tabaghi2021, Kriouko2010, Moshir2021, Wilson2014}.

The advancement of optimization techniques in Euclidean space, particularly for handling
constrained problems, encompasses a wide spectrum of methodologies, each specifically
tailored to address different types of constraints. A thorough understanding of each
constraint type is pivotal in determining the most suitable methodology for its resolution.
For example, a common methodology employed in solving non-linear programming problems is the
Augmented Lagrangian or Sequential Quadratic Methods. However, these methods require a deep
study of Lagrange multipliers, generate sequences that are infeasible, and only achieve the
desired feasibility in the limit. It is worth noting that when constraints are considered
easy, meaning that projection calculation can be efficiently done or has a formula,
penalizing these constraints may not be an advantageous strategy, as explained in
\cite{Andreanietall2007}.  Another employed methodology is the Frank-Wolfe strategy, which
involves minimizing the linearization of the objective function over a compact constraint set.
This type of methodology requires an oracle that minimizes a linear function over compact and
easy-to-handle constraints. A very popular methodology is also interior point methods. These
methods generally parameterize optimality conditions to define subproblems, and for each
parameter, they use the Newton method to find a solution of it. However, obtaining solutions
for the subproblems is very costly even for very simple constraint sets, due to the
application of the Newton method iteration. An alternative to all these methods, which can be
advantageous in many aspects, is the use of projected-type methods, such as the gradient
projection method, when the calculation of the projection  is straightforward or has a
formula. In this case, the sequence produced by the method is feasible,  there is no need to
deal with multipliers and Newton's equations, no requirement for an oracle, and furthermore,
the constraint sets do not need to be compact.   Optimization problems constrained to a
Riemannian manifold are typically reformulated as unconstrained problems by treating the
constraint set as the ambient space. Consequently, numerous unconstrained optimization
methods developed for Euclidean spaces have been adapted to address this novel context.  For
instance, methods such as the gradient method \cite{wilson2018gradient} and the Riemannian trust-region algorithm \cite{jawanpuria2019low} have been tailored for application in  hyperbolic spaces. Additionally, there are many more examples across various Riemannian
manifolds, as seen in books such as \cite{Absil2008, Boumal2020}. However, when additional
constraints are imposed directly within the manifold, the problems become intrinsically
constrained. In such scenarios, optimization problems on manifolds with intrinsic constraints
pose significantly greater challenges.   In theory, similar to how it was done for
constrained optimization, we can adapt methodologies from Euclidean space to address
intrinsically constrained problem in Riemannian manifolds.  Various methods have already been
developed to address intrinsically constrained problems in Riemannian manifolds.  These
encompass the Augmented Lagrangian method \cite{LiuBoumal2020}, Sequential Quadratic Methods
\cite{Obara2022}, Frank-Wolfe algorithm \cite{WeberSra2023}, and interior point method
\cite{LaiYoshise2024}. While these adaptations have successfully overcome challenges inherent
to constrained optimization, it is crucial to recognize that they may still inherit certain
deficiencies from their Euclidean counterparts, as discussed above. Additionally, new
challenges may arise in this new scenario. For example, in the Frank-Wolfe algorithm, the
subproblem in Riemannian settings is not necessarily linear or convex. Consequently, it
generally shares the same computational complexity as the original problem, making it
applicable efficiently only to very specific types of problems. Notably, discussions
comparing the advantages of different methodologies in Euclidean spaces remain pertinent in
this context. It is important to emphasize that the theoretical foundations required to
support and design methods for intrinsic constrained optimization on Riemannian settings are
still in the early stages of development. For instance, a fundamental challenge encountered
when adapting projected-type methods is computing the intrinsic projection into geodesically
convex sets
of a Riemannian manifold, a problem that remains systematically unresolved even for simple
geodesically convex sets. The problem of computing the intrinsic projection into geodesically
convex sets was initially
addressed in the paper \cite{Walter1974}, where intriguing theoretical properties were
discovered. Due to the difficulty of addressing the problem in a general manner, the idea is
to concentrate efforts on computing the intrinsic projection into specific geodesically convex sets of a
Riemannian manifold. For example, papers  \cite{Bauschke2018, FerreiraIusemSandor2013}  have
studied the problem of projecting onto an intrinsic spherically convex subset of the sphere,
where an
explicit formula for the intrinsic projection was reduced to the Euclidean projection on a
suitable cone. We anticipate that this strategy will become a trend in the near future.  In
our current paper, we adopt a similar strategy to tackle the problem of intrinsic projection
into hyperbolically convex sets of hyperbolic space forms.

As previously discussed,  efficiently computing the intrinsic projection into geodesically
convex sets or
deriving formulas for specific geodesically convex sets of a Riemannian manifold remains a challenging
task. Consequently, these computational obstacles severely limit the practical exploration of
numerical projected-type methods aimed at addressing intrinsic constrained problems within
Riemannian manifolds. For example, successfully computing the intrinsic projection into
geodesically convex sets would facilitate the practical application of the Alternated
Projection Methods investigated in \cite{Lewis2022}, where only theoretical properties are presented. Similarly, this limitation applies to the Gradient Projection Method investigated in \cite{zhang2016},
where the presented numerical results are solely for an unconstrained problem.  Currently,
there are only a few studies addressing the gradient projection method in the Riemannian
context, despite its significance as one of the foremost first-order methods for constrained
optimization. Previous research has explored the application of the gradient projection for
minimizing functions on Riemannian manifolds under convexity or coercivity assumptions, as
well as the boundedness of the set of constraints. However, even for this method, which plays
a fundamental role in constrained optimization, the conducted studies predominantly focus on
theoretical findings regarding the computational complexity associated with solving problems
with convex or coercive objective functions, as seen in
\cite{criscitiello2022negative,Feng2022,HuangWei2022,Martinez023, zhang2016}. This focus is
understandable due to the limitation of computing the intrinsic projection into geodesically
convex sets.

This paper, a continuation of \cite[Section~3]{FerreiraNemethShu2022}, addresses the natural counterpart in hyperbolic space forms of the problem
presented in \cite{FerreiraIusemSandor2013}, which involves computing a projection formula into a spherically convex set of the sphere. However,
as we will see, the problem in hyperbolic space forms is considerably more difficult. This work is the first to systematically derive and analyze explicit formulas for projection mappings specifically within the context of hyperbolic space forms. Additionally, we propose and analyze two intrinsic versions
of the gradient projection method, one with constant step sizes and the other with backtracking step sizes, to solve problems constrained to
$\kappa$-hyperbolically convex sets and for potentially non-convex objective functions in $\kappa$-hyperbolic space forms with a constant
sectional curvature $-\kappa<0$.  For our study, we selected the hyperboloid model (also known as the Lorentz model) from among the existing
hyperbolic geometry models (see \cite{BenedettiPetronio1992, Cannon1997}), as it appears to be better suited for studying Riemannian optimization
problems and provides more numerical stability to solution methods (see \cite{nickel2018learning, peng2021hyperbolic}).  We begin by reviewing
some fundamental properties of $\kappa$-hyperbolic space forms, followed by presenting novel properties of the intrinsic $\kappa$-projection into
$\kappa$-hyperbolically convex sets of $\kappa$-hyperbolic space forms.  These properties are crucial for analyzing the gradient projection method
and hold intrinsic importance on their own. For instance, thanks to these properties, we can conduct the analysis of the gradient projection
method without the need to assume that the constraint set is bounded. Specifically, we introduce the concept of the Lorentz projection and
establish the relationship between the intrinsic $\kappa$-projection, the Lorentz projection, and the Euclidean orthogonal projection.
Additionally, we provide formulas for the intrinsic $\kappa$-projection into specific $\kappa$-hyperbolically convex sets using the Euclidean orthogonal projection and the Lorentz projection. It is worth noting that these formulas enable practical applications to problems with these sets as constraints. Concerning the convergence results obtained from the gradient projection method, we demonstrate that every accumulation point of the sequence generated by the method with backtracking step sizes is a stationary point for the problem. Under the assumption of Lipschitz continuity of the gradient of the objective function, we show that each accumulation point of the sequence generated by the gradient projection method with constant step size is a stationary point. Furthermore, we provide an iteration complexity bound that characterizes the number of iterations needed to achieve a suitable measure of stationarity for both step sizes.   It is important to emphasize that the findings presented here do not rely on the convexity or coercivity of the objective function, nor on the compactness of the constraint set,  as previously considered, for example, in  \cite{zhang2016}. Finally, we explore the properties of the constrained version of the problem of computing the center of mass, commonly referred to as the Fermat-Weber problem. Specifically, we demonstrate that the sequence generated by the gradient projection method with backtracking step sizes converges to the unique solution of this problem. To complement the theoretical developments, we present numerical experiments in Section~\ref{sec:Nuerics}, conducted to solve the Fermat-Weber problem. These experiments validate our theoretical findings and illustrate the practical performance of the proposed gradient projection methods in hyperbolic space forms.

The structure of this paper is as follows: In Section~\ref{sec:int.01}, we review some notations and basic results. Section~\ref{sec:int.1}
provides an overview of the geometry of ${\kappa}$-hyperbolic space forms, including notations, definitions, and basic results.
Section~\ref{sec:pj} presents some properties of the intrinsic $\kappa$-projection into $\kappa$-hyperbolically convex sets, with a focus on
explicit computation of the intrinsic $\kappa$-projection for specific $\kappa$-hyperbolically convex sets.
In Section~\ref{sec:socop}, we address the constrained optimization problem and introduce two versions of the gradient projection method to solve it.
Section~\ref{sec:App} explores  some  properties of the constrained version of the problem of computing the center of mass.
Section~\ref{sec:Nuerics} explores different algorithms for this case numerically. Finally, Section~\ref{sec:Conclusions} concludes the paper.
\subsection{Notation} \label{sec:int.01}
 Let ${\mathbb R}^{m}$ be   the $m$-dimensional Euclidean space.  The set of all
 $m \times n$ matrices with real entries is denoted by ${\mathbb R}^{m \times n}$
 and ${\mathbb R}^m\equiv {\mathbb R}^{m\times 1}$. The canonical vectors  of
 ${\mathbb R}^m$ are denoted by  $e^1,\dots,e^m$.  For   $M \in {\mathbb
 R}^{m\times n}$ the matrix $M^{\top}  \in {\mathbb R}^{n\times m}$  denotes  the
 {\it transpose} of $M$.     The matrix ${\rm I}$ denotes the $n\times n$
 identity matrix.  Denote
 $\R^m_+=\{x=(x^1,\dots,x^m)^T\in\R^m:x_1\ge0,\dots,x^m\ge0\}$ the nonnegative
 orthant of ${\mathbb R}^m$. Let $\{v^1, v^2, \ldots, v^m\}$ be a set of vectors of  $\R^m$. Then, the cone defined by
$
\cone\left\{v^1,\dots,v^m\right\}:=\left\{\lambda_1v^1+\cdots+\lambda_mv^m: ~\lambda_1\ge0,\dots,\lambda_m\ge0\right\},
$
 is called a  \emph{simplicial cone} or \emph{finitely generated cone}.   Let ${K} \subseteq {\mathbb R}^{n+1}$  be a  convex cone. The {\it Euclidean orthogonal projection set-valued mapping    $\Pi_{K}:{\mathbb R}^{n+1} \multimap {K}$ into  ${K}$} is defined by
\begin{equation}\label{eq:pirn}
\Pi_{K}(x):=\Big\{y\in {\mathbb R}^{n+1}:~y=\underset{z\in {K}}{\argmin}\sqrt{(x-z)\tp(x-z)}\Big\}.
\end{equation}
Since   ${K}$ is convex, $\Pi_{K}(x)$ is either an empty set or contains a single point. Moreover,  it is well known that, if \(y\in  \Pi_{K}(x)\) then
\begin{equation}\label{eq:zara}
	(x-y)^{^\top}z\leq 0,\mbox{ }  \forall\, z\in {K}\qquad\mbox{and}\qquad (x-y)^\top y=0,
\end{equation}
see for example \cite[Theorem 3.1.1]{HiriartLemarechal2001}.  Whenever  the set
${K}$ is  closed and convex, the projection set is nonempty and  has a unique point, which is
characterised by \eqref{eq:zara}. In this case, we use the notation
$
\Pi_{K}(x)={\argmin}_{{z\in {K}}}\sqrt{(x-z)\tp (x-z)}.
$
Let $\alpha\ge 1$ be  a real constant.  Denote de {\it $\alpha$-circular cone} ${L}_{\alpha}$ by
\begin{equation}\label{eq:cc}
		{L}_{\alpha}=\Big\{(x_1,\dots,x_n,x_{n+1})\in\R^{n+1}:~x_{n+1}\ge\alpha\sqrt{x_1^2+\dots+x_n^2}\Big\}.
\end{equation}
 When $\alpha=1$, the  $\alpha$-circular cone ${L}_{\alpha}$ is known as the {\it Lorentz cone} and will simply be denoted by ${L}$.
 Note that, for all $\alpha>1$, we have  ${L}_{\alpha}\subseteq \inte(\Lo)$ and
 it is a closed and convex cone.
\section{Basics results about the $\kappa$-hyperbolic space forms} \label{sec:int.1}
In this section, we provide  notations, definitions, and basic results related to the geometry of the ${\kappa}$-hyperbolic space forms that are utilized throughout the paper. This section constitutes a review of the helpful  ${\kappa}$-hyperbolic space forms results used in optimization. The books  that we refer to in this section are \cite{BenedettiPetronio1992} and \cite{Ratcliffe2019}, and we also draw insights from \cite{Boumal2020}.

Let $\lng\cdot , \cdot \rng$ be the {\it Lorentzian inner product}  of  $x:=(x_1, \ldots,x_n, x_{n+1})^{\top} $ and  $y:=(y_1, \ldots, y_n,y_{n+1})^{\top}$ on ${{\mathbb R}^{n+1}}$ defined  as follows
\begin{equation} \label{eq:ip}
\lng x, y\rng :=x^{\top}{\rm J}y = x_{1}y_{1}+ \cdots + x_{n}y_{n}-x_{n+1}y_{n+1},
\end{equation}
where  the diagonal matrix  ${\rm J}$ is defined by ${\rm J}:={\rm diag}(1, \ldots,1, -1) \in {\mathbb R}^{(n+1)\times (n+1)}$. For each  $x\in {{\mathbb R}^{n+1}}$,  the {\it Lorentzian norm (length)} of $x$ is  the complex number $\|x\|:= \sqrt{\lng x, x\rng}.$
\begin{remark} \label{re;PLC}
From \eqref{eq:cc}  with $\alpha=1$ and  \eqref{eq:ip} we have ${L}:=\{x:=(x_1, \ldots, x_{n+1})^T\in  {\mathbb R}^{n+1}:~ \lng x, x\rng\leq 0,~ x_{n+1}\geq 0\}.$
Furthermore, $ \lng x, y\rng< 0$,  for all $x\in\Lo$ and all $y\in\inte{L}=\{x:=(x_1, \ldots, x_{n+1})^T\in  {\mathbb R}^{n+1}:~ \lng x, x\rng< 0,~ x_{n+1}>0\},$
because $\lng x,y\rng=-y\tp(-\J x)$, $-\J x\in -\J\Lo=\Lo$ and $\Lo$ is
self-dual.
\end{remark}
The {\it $n$-dimensional  ${\kappa}$-hyperbolic space form} and its {\it tangent hyperplane at  $p$} are denoted  by
\begin{equation} \label{eq:hs}
{\mathbb H}^{n}_{\kappa}:=\Big\{ p\in {{\mathbb R}^{n+1}}:~\langle p, p\rangle
=-\frac{1}{\kappa}, ~p^{n+1}>0\Big\}, \qquad \quad T_{p}{{\mathbb H}^n_{\kappa}}:=\Big\{v\in {{\mathbb R}^{n+1}}\, :\,
\lng p, v \rng=0\Big\},
\end{equation}
respectively, for all $\kappa>0$.
The  Lorentzian inner product   is not  positive definite in the entire space ${{\mathbb R}^{n+1}}$.  However,
one can show that its restriction to the tangent spaces of ${\mathbb H}^{n}_{\kappa}$ is positive
definite; see \cite[Section 7.6]{Boumal2020}.  The  {\it  Lorentzian projection into} $T_p{{\mathbb H}^n_{\kappa}}$ is the linear mapping $\Proj^{\kappa}_p:{{\mathbb R}^{n+1}} \to T_p{{\mathbb H}^n_{\kappa}}$  defined by
\begin{equation} \label{eq: proj}
 \Proj^{\kappa}_p x:=x+ \kappa \langle p, x \rangle p,
\end{equation}
 i.e., $ \Proj^{\kappa}_p:={\rm I} +\kappa pp^{\top}{\rm J} $, where ${\rm I}  \in {\mathbb R}^{(n+1)\times (n+1)}$ is the identity matrix.
\begin{remark}
The  Lorentzian projection \eqref{eq: proj} is self-adjoint with respect to the   Lorentzian inner product \eqref{eq:ip}. Indeed, $\lng \Proj^{\kappa}_pu, v\rng=\lng u, \Proj^{\kappa}_pv\rng$,  for all   $u, v \in  {{\mathbb R}^{n+1}}$ and all $p\in  {{\mathbb H}^n_{\kappa}}$. Moreover, we also have $\Proj^{\kappa}_p\Proj^{\kappa}_p=\Proj^{\kappa}_p$, for all $p\in  {{\mathbb H}^n_{\kappa}}$.
\end{remark}
The {\it intrinsic $\kappa$-distance on the  $\kappa$-hyperbolic space form} between two  points $p, q \in {{\mathbb H}^n_{\kappa}}$  is  given~by
\begin{equation} \label{eq:Intdist}
d_{\kappa}(p, q):=\frac{1}{\sqrt{\kappa}}{\rm arcosh} (-\kappa \lng p , q\rng).
\end{equation}
 It can  be shown that   $({{\mathbb H}^n_{\kappa}}, d_{\kappa})$  is a complete metric space, so that $d_{\kappa}(p,q)\ge 0$ for all $p,q \in {{\mathbb H}^n_{\kappa}}$, and
$d_{\kappa}(p,q)=0$ if and only if $p=q$.   Moreover,  $({{\mathbb H}^n_{\kappa}}, d)$ has the same topology as ${{\mathbb R}^{n}}$ and $({{\mathbb H}^n_{\kappa}}, \lng\cdot , \cdot \rng)$ is a Riemannian manifold with constant sectional   curvature $-\kappa<0$.  Next we give two standard notations. We denote the {\it open} and {\it closed ball} of radius $\delta >0$ and center
$p\in {{\mathbb H}^n_{\kappa}}$ by $B_{\delta}(p):=\{q\in {{\mathbb H}^n_{\kappa}} : d_{\kappa}(p,q)<\delta
\}$  and $\bar{B}_{\delta}(p):=\{q\in {{\mathbb H}^n_{\kappa}} : d_{\kappa}(p,q)\leq \delta \}$, respectively. The intersection curve of a plane though the origin of ${{\mathbb R}^{n+1}}$
with  $ {{\mathbb H}^n_{\kappa}}$ is called a { \it geodesic}. If $ p, q\in {{\mathbb H}^n_{\kappa}}$ and  $q\neq p$, then the   unique {\it geodesic segment  from $p$ to $q$ } is
\begin{equation*}
\gamma_{pq}(t)= \Big( \cosh(t\sqrt{\kappa}) + \frac{{\kappa}\lng p, q\rng \sinh(t\sqrt{\kappa}) }{\sqrt{{\kappa}^2\lng
p, q\rng^2-1}}\Big) p
+ \frac{\sinh(t\sqrt{\kappa}) }{\sqrt{{\kappa}^2\lng p, q\rng^2-1}}\;q, \qquad \forall t\in [0, \;d_{\kappa}(p,q)].
\end{equation*}
The {\it $\kappa$-exponential mapping} $\exp^{\kappa}_{p}:T_{p}{{\mathbb H}^n_{\kappa}} \rightarrow {{\mathbb H}^n_{\kappa}} $ is defined
by $\exp^{\kappa}_{p}v= \gamma _{v}(1)$, where $\gamma _{v}$ is the geodesic defined by its
initial position $p$, with velocity $v$ at $p$. Hence, $\exp^{\kappa}_{p}v=p$ for  $v=0$.   It is easy to prove that $\gamma _{tv}(1)=\gamma _{v}(t)$ for any values of $t$.
Therefore,  for all $t\in {\mathbb R}$ we have
\begin{equation} \label{eq:geoexp}
\exp^{\kappa}_{p}(tv):=  \displaystyle \cosh(t\sqrt{\kappa}\|v\|) \,p+ \sinh(t\sqrt{\kappa}\|v\|)\, \frac{v}{\sqrt{\kappa}\|v\|},  \qquad  \forall v\in T_p{{\mathbb H}^n_{\kappa}}/\{0\}.
\end{equation}
We will also use the expression above for denoting  the geodesic starting
at $p\in {{\mathbb H}^n_{\kappa}}$ with velocity $v\in T_p{{\mathbb H}^n_{\kappa}}$ at $p$.
The {\it inverse of the $\kappa$-exponential mapping}  is given by $\log^{\kappa}_{p}q=0$, for $q=p$, and
\begin{equation} \label{eq:expinv}
\log^{\kappa}_{p}q:=  \displaystyle \frac{\sqrt{\kappa}d_{\kappa}(p, q)}{\sqrt{{\kappa}^2\lng p, q\rng^2-1}}
\Proj^{\kappa}_pq= \displaystyle  d_{\kappa}(p, q)\frac{\Proj^{\kappa}_pq}{\|\Proj^{\kappa}_pq\|},  \qquad  q\neq p.
\end{equation}
It follows from \eqref{eq:expinv} that
\begin{equation} \label{eq:edn}
d_{\kappa}(p, q)=\|\log^{\kappa}_{p}q\|, \qquad  p, q \in {{\mathbb H}^n_{\kappa}}.
\end{equation}
In the next lemma, we recall the well-known ``cosine law'' for triangles in the $\kappa$-hyperbolic space form.
\begin{lemma}  \label{le:CosLaw1}
  Let   ${x}, {y}, {z} \in {{\mathbb H}^n_{\kappa}}$. Let $\theta_{x}$  be the angle between the vectors $ \log^{\kappa}_{{x}}{y}$ and $\log^{\kappa}_{x}{z}$. Then,
  \begin{equation*}
    \cosh(\sqrt{\kappa} d_{\kappa}({y}, {z}))=\cosh(\sqrt{\kappa} d_{\kappa}({x}, {y}))\cosh(\sqrt{\kappa} d_{\kappa}({x}, {z}))-\sinh(\sqrt{\kappa} d_{\kappa}({x}, {y})) \sinh(\sqrt{\kappa} d_{\kappa}({x}, {z}))  \cos \theta_{x}.
  \end{equation*}
 \end{lemma}
Next, we will recall the well-known ``comparison theorem" for triangles in the $\kappa$-hyperbolic space form, as stated in \cite[Proposition 4.5]{Sakai1996}. Its proof is a straightforward application of Lemma~\ref{le:CosLaw1}.
 \begin{lemma}  \label{le:CosLaw2}
There   holds $d_{\kappa}^2({x}, {y})+d_{\kappa}^2({x},{z})-2\lng  \log^{\kappa}_{{x}}{y}, \log^{\kappa}_{{x}}{z}\rng\leq d_{\kappa}^2({y},{z})$,  for all  ${x}, {y}, {z} \in {{\mathbb H}^n_{\kappa}}$.
\end{lemma}
Let $\Omega \subseteq {{\mathbb H}^n_{\kappa}}$ be an  open set and $f: \Omega
\to  {\mathbb R}$ a differentiable function.  The {\it gradient on the
${\kappa}$-hyperbolic space form} of $f$  is the unique vector field $\Omega \ni p\mapsto \grad f(p)\in T_{p}{M}$ such that $df(p)v=\lng \grad f(p), v\rng$, see \cite[Proposition~7-5, p.162]{Boumal2020}.  Therefore, we have
\begin{equation} \label{eq:grad}
\grad f(p):= \Proj^{\kappa}_p ({\rm J}f'(p))=  {\rm J}f'(p)+ \kappa \lng {\rm J}f'(p), p\rng \,p,
\end{equation}
where $f'(p) \in {\mathbb R}^{n+1}$ is the usual gradient of  $ f$ at $p\in \Omega$.    The {\it Hessian on the ${\kappa}$-hyperbolic space form}
of a twice differentiable function $f: \Omega \to  {\mathbb R}$ at $p\in \Omega$ is  the  mapping
$ \Hess f(p):T_p{{\mathbb H}^n_{\kappa}} \to T_p{{\mathbb H}^n_{\kappa}}$ given by
\begin{equation} \label{eq:Hess}
\Hess f(p):= \Proj^{\kappa}_p\left({\rm J} f''(p)+{\kappa}\lng {\rm J}f'(p), p\rng {\rm I}\right),
\end{equation}
where $ f''(p)$ is the usual Hessian of  $f$ at  $p$, see  \cite[Proposition~7.6, p.163]{Boumal2020}.  For each $p, q\in {{\mathbb H}^n_{\kappa}}$, the {\it parallel transport} mapping  $P_{pq} \colon T _{p} {{\mathbb H}^n_{\kappa}} \to T _ {q}{{\mathbb H}^n_{\kappa}}$  is given by
\begin{equation*}
{P}_{pq}(v):= v-\frac{\lng v,  \log^{\kappa}_{q}p \rng}{d_{\kappa}^2(p, q)}\left( \log^{\kappa}_{q}p+ \log^{\kappa}_{p}q \right)=\left[{\rm I}-\frac{\kappa}{{\kappa}\lng p, q\rng-1}(p+q)q^\top {\rm J}\right]v.
\end{equation*}
Note that  for any geodesic segment  $\gamma:[a,b] \to  {{\mathbb H}^n_{\kappa}}$ we have
$ \gamma'(t)={P}_{\gamma(a)\gamma(t)}(\gamma'(a))$, for all $t\in [a, b]$, or equivalently
$\gamma''(t)=0$,  for all $t\in [a, b]$.
\begin{definition} \label{Def:GradLips}
 Let   \( f: {D} \to {\mathbb R}\)  be a differentiable function and  ${D}\subseteq {{\mathbb H}^n_{\kappa}}$. The gradient vector  field  $\grad f$   is  said  to be   Lipschitz continuous on ${D}$ with   constant $L\geq 0$ if  the following inequality holds $\left\|P_{pq} \grad f(p)- \grad f(q)\right\|\leq Ld_{\kappa}(p,q)$,   for all  $p,q\in{D}$.
\end{definition}
The next lemma can be found,  for example,  in  \cite[Corollary 10.47]{Boumal2020}.
\begin{lemma} \label{le:CharactGL}
 Let   \( f: {D} \to {\mathbb R}\)  be  a twice continuously differentiable function and  ${D}\subseteq {{\mathbb H}^n_{\kappa}}$. The gradient vector  field  $\grad f$   is  Lipschitz continuous on ${D}$ with   constant $L\geq 0$ if  and only if  there exists $L\geq 0$ such that $\|\mbox{Hess}\,f(p)\|\le L$, for all $p\in {D}$.
\end{lemma}
The set ${C}\subseteq{{\mathbb H}^n_{\kappa}} $ is called  \emph{$\kappa$-hyperbolically convex} if for any $p$, $q\in {C}$   the geodesic segment joining
$p$ to $q$  is contained in ${C}$. For each set $ {A} \subseteq {{\mathbb H}^n_{\kappa}}$, let ${K_A}$ be the {\it cone spanned by} ${A}$,
namely,
\begin{equation} 
{K_A}:=\left\{ tp \, :\, p\in A, \; t\in [0, +\infty) \right\}.
\end{equation}
Clearly, ${K_A}$ is the smallest cone which contains ${A}$ and ${K_A}\setminus\{0\}$ belongs to the interior of the Lorentz cone ${L}$.
In the next proposition  we relate a $\kappa$-hyperbolically convex set with the cone spanned by it, its proof can be found in  \cite[Proposition 1]{FerreiraNemethShu2022}.

Subject to a few minor adjustments, the next lemma closely mirrors  \cite[Proposition 10.54]{Boumal2020}

\begin{lemma} \label{le:lc}
Let ${D} \subseteq {{\mathbb H}^n_{\kappa}}$ and $f: {D} \to {\mathbb R}$ be a differentiable function. If ${D}$ is a subset of ${\mathbb
H}^n_{\kappa}$ and $\grad f$ is Lipschitz continuous on a $\kappa$-hyperbolically convex set ${C} \subseteq {D}$ with a constant $L \geq 0$, then for all $p, q \in {C}$, the following inequality holds
$f( q) \leq f(p) + \lng\grad f(p), \log^{\kappa}_{p}q\rng +({L}/{2})d_{\kappa}^2(p,q)$, for all $p, q\in{C}$.
\end{lemma}

\begin{proposition} \label{pr:ccs}
The set ${C}\subseteq\mathbb H^n$ is $\kappa$-hyperbolically convex if and only if the cone
${{K}_{C}}$ is convex.
\end{proposition}
\begin{remark}
	The $\kappa$-hyperbolically closed  convex sets are intersections of the hyperboloid  with  convex cones which belong to the Lorentz cone  ${L}$.
Indeed, it follows  from Proposition~\ref{pr:ccs}, that if ${K}\subseteq {L}$ is a  convex cone, where ${L}$ is the  Lorentz
cone,  then ${C}={K}\cap {{\mathbb H}^n_{\kappa}}$ is a $\kappa$-hyperbolically  convex set and ${K}={{K}_{C}}$.
\end{remark}

\section{Intrinsic $\kappa$-projection into $\kappa$-hyperbolically convex sets}\label{sec:pj}
In this section we present new  properties of the  intrinsic $\kappa$-projection into $\kappa$-hyperbolically convex sets on ${\kappa}$-hyperbolic space forms. Let ${{C}\subseteq} {{\mathbb H}^n_{\kappa}} $ be a closed $\kappa$-hyperbolically convex set and $p\in{{\mathbb H}^n_{\kappa}}$.  Consider the following constrained  optimization problem
\begin{align} \label{d:copf}
 \min_{q\in {C}}  d_{\kappa}(p,q).
\end{align}
The minimal value of the function ${C}  \ni q \mapsto d_{\kappa}(p,q)$ is called the
{\it distance of $p$ from  ${C}$} and it is denoted by $ d_{{\kappa},{C}} $,
i.e.,  $ d_{{\kappa},{C}} : {{\mathbb H}^n_{\kappa}} \to  {\mathbb R} $ is defined by
$
 d_{{\kappa},{C}}(p):=\min_{q\in {C}}  d_{\kappa}(p,q).
$
Since $({{\mathbb H}^n_{\kappa}}, d_{\kappa})$  is a complete metric space,  \cite[Proposition 2]{FerreiraNemethShu2022} implies that   $ | d_{{\kappa},{C}}(p) -  d_{{\kappa},{C}}(q)|\leq d_{\kappa}(p, q)$,  for all  $p, q \in {{\mathbb H}^n_{\kappa}}$. In particular, the function $ d_{{\kappa},{C}}$ is continuous.  Note that the problem \eqref{d:copf} has a solution because ${C}$ is a closed set and the distance function is continuous.
The solution of problem \eqref{d:copf} is referred to as the ``metric $\kappa$-projection",  which was first studied in \cite{Walter1974}. Similar to the proof in \cite[Proposition~4.2]{FerreiraNemethShu2022}, we can also demonstrate that the ``intrinsic $\kappa$-projection mapping" ${P}^{\kappa}_{C}: {{\mathbb H}^n_{\kappa}}  \to {C}$ defined by
\begin{align} \label{d:prjection}
	{P}^{\kappa}_{C}(p):= \underset{q\in {C}}{\argmin}\,d_{\kappa}(p,q).
\end{align}
is well-defined. In other words, for each $p\in {{\mathbb H}^n}$, the minimizer ${P}^{\kappa}_{C}(p)$ of problem \eqref{d:prjection} exists and is unique. The following proposition can also be proven in a similar manner to \cite[Proposition 4.4]{FerreiraNemethShu2022}.
\begin{proposition} \label{pr:PrjCont}
Let ${{C}\subseteq} {{\mathbb H}^n}$ be a nonempty closed $\kappa$-hyperbolically convex set. Then ${P}^{\kappa}_{C}$ is continuous.
\end{proposition}

Note that the definition provided for the intrinsic $\kappa$-projection mapping \eqref{d:prjection} implies the inequality
\begin{equation}\label{d:prj}
d_{\kappa}(p,{P}^{\kappa}_{C}(p))\leq d_{\kappa}(p,q),
\qquad   \forall\, q\in {C}.
\end{equation}
Furthermore, we can prove that equations \eqref{d:prjection} and \eqref{d:prj} are equivalent to the following inequality
\begin{equation} \label{d:prjectionef}
\langle \Proj^{\kappa}_{{P}^{\kappa}_{C}(p)}p, \Proj^{\kappa}_{{P}^{\kappa}_{C}(p)}q \rangle \leq 0, \quad \forall q\in {C}, \forall p\in {{\mathbb H}^n_{\kappa}}.
\end{equation}
By considering that, for all $p, q \in {{\mathbb H}^n_{\kappa}}$, we have $\langle {P}^{\kappa}_{C}(p), q \rangle \leq -1/{\kappa}$, we can deduce from \eqref{eq:expinv} that \eqref{d:prjectionef}   can alternatively be expressed as follows
\begin{equation} \label{d:uexp}
\langle \log^{\kappa}_{{P}^{\kappa}_{C}(p)}p, \log^{\kappa}_{{P}^{\kappa}_{C}(p)}q \rangle \leq 0, \quad \forall q\in {C}, \forall p\in {{\mathbb H}^n_{\kappa}},
\end{equation}
as shown in \cite[Corollary 3.1]{FerreiraOliveira2002}. Furthermore, since the function $[1, +\infty] \ni \tau \mapsto {\rm arcosh}(\tau)$ is increasing, it follows from \eqref{eq:Intdist} that \eqref{d:prjection}, \eqref{d:prjectionef}, and \eqref{d:uexp} are also equivalent to
\begin{align} \label{d:prjcef}
	{P}^{\kappa}_{C}(p):= \underset{q\in {C}}{\argmin}\, (-\langle p, q \rangle).
\end{align}
Therefore, by combining \eqref{eq:Intdist} with \eqref{d:prj}, we conclude that \eqref{d:prjcef} can be stated equivalently as follows
\begin{equation} \label{d:prj2}
\langle p, q \rangle \leq \langle p, {P}^{\kappa}_{C}(p) \rangle, \quad \forall q\in {C}.
\end{equation}
As an application of Lemma~\ref{le:CosLaw2}, we have the following property of the intrinsic $\kappa$ projection.
\begin{lemma} \label{le:coslawap}
Let ${C} \subseteq{{\mathbb H}^n_{\kappa}}$ be a closed $\kappa$-hyperbolically convex set, and let ${q}\in{{\mathbb H}^n_{\kappa}}$. Then, there holds
\begin{equation}\label{eq;dlcf0}
d_{\kappa}^2(p, {P}^{\kappa}_ {C}({q}))\leq \langle \log^{\kappa}_{p}{q}, \log^{\kappa}_{p} {P}^{\kappa}_ {C}({q})\rangle
+\langle \log^{\kappa}_{{P}^{\kappa}_{C}(q)}p, \log^{\kappa}_{{P}^{\kappa}_{C}(q)}q\rangle , \quad \forall p\in {C}.
\end{equation}
As a consequence, we have
\begin{equation}\label{eq;dlcf}
d_{\kappa}^2(p, {P}^{\kappa}_ {C}({q}))\leq \langle \log^{\kappa}_{p}{q}, \log^{\kappa}_{p} {P}^{\kappa}_ {C}({q})\rangle, \quad \forall p\in {C}.
\end{equation}
Furthermore, we have $d_{\kappa}(p, {P}^{\kappa}_ {C}({q}))\leq d_{\kappa}(p, q)$ for all $p\in {C}$.
\end{lemma}
\begin{proof}
Let $p\in {C}$ and $q\in{{\mathbb H}^n_{\kappa}}$. By applying Lemma~\ref{le:CosLaw2} with $x=p$, $y={P}^{\kappa}_ {C}(q)$ and $z=q$, we obtain:
\begin{equation} \label{eq:coslaw1}
d_{\kappa}^2(p, {P}^{\kappa}_ {C}(q))+d_{\kappa}^2(p,q)-2\langle \log^{\kappa}_{p}{P}^{\kappa}_ {C}(q), \log^{\kappa}_{p}q\rangle \leq d_{\kappa}^2({P}^{\kappa}_ {C}(q),q).
\end{equation}
Similarly, by applying Lemma~\ref{le:CosLaw2} with $x={P}^{\kappa}_ {C}(q)$, $y=p$ and $z=q$, we have:
\begin{equation} \label{eq:coslaw2}
d_{\kappa}^2({P}^{\kappa}_ {C}(q), p)+d_{\kappa}^2({P}^{\kappa}_
{C}(q),q)-2\langle \log^{\kappa}_{{P}^{\kappa}_{C}(q)}p, \log^{\kappa}_{{P}^{\kappa}_{C}(q)}q\rangle \leq d_{\kappa}^2(p,q).
\end{equation}
By combining \eqref{eq:coslaw1} with \eqref{eq:coslaw2}, we conclude that
\eqref{eq;dlcf0} holds. Therefore, since $p\in {C}$, by applying \eqref{d:uexp}, the inequality \eqref{eq;dlcf} follows from \eqref{eq;dlcf0}. The last inequality of the lemma is a direct consequence of the second Cauchy inequality and \eqref{eq:edn}.
\end{proof}
Based on Lemma~\ref{le:coslawap}, we next derive a crucial property of the intrinsic $\kappa$-projection that holds significant importance in the analysis of the gradient projection method.
\begin{proposition}  \label{pr:condw2}
 Let   $p \in {C}$, $0\neq v\in T_{p}{{\mathbb H}^n_{\kappa}}$ and $\alpha >0$.  Then, we have
 \begin{equation} \label{eq;snlpf}
 \left\langle {v},  \log^{\kappa}_{p}{P}^{\kappa}_{C}\left(\exp^{\kappa}_{p}\left(-\alpha {v}\right)\right) \right\rangle \leq - \frac{1}{\alpha}d_{\kappa}^2\left(p,{P}^{\kappa}_{C}(\exp^{\kappa}_{p}(-\alpha {v}))\right).
\end{equation}
\end{proposition}
\begin{proof}
For simplicity, let us define $\gamma(\alpha):=\mbox{exp}_{p}(-\alpha {v})$. Thus,  $v=-({1}/{\alpha})\log^{\kappa}_{p}\gamma(\alpha)$. After performing some algebraic manipulations, we can obtain the following expression
\begin{equation*}
\left\langle {v}, \log^{\kappa}_{p}{P}^{\kappa}_{C}(\gamma(\alpha)) \right\rangle =
-\frac{1}{\alpha}d_{\kappa}^2 (p, {P}^{\kappa}_{C}(\gamma(\alpha))) +
\frac{1}{\alpha}\left(d_{\kappa}^2 (p, {P}^{\kappa}_{C}(\gamma(\alpha)))- \left\langle \log^{\kappa}_{p}\gamma(\alpha), \log^{\kappa}_{p}{P}^{\kappa}_{C}(\gamma(\alpha))\right\rangle\right).
\end{equation*}
Since $p \in {C}$,   Lemma~\ref{le:coslawap} with $q=\gamma(\alpha)$ implies $d_{\kappa}^2 (p, {P}^{\kappa}_{C}(\gamma(\alpha))) - \langle \log^{\kappa}_{p}\gamma(\alpha), \log^{\kappa}_{p}{P}^{\kappa}_{C}(\gamma(\alpha))\rangle\leq 0$.  Thus,   it follows from  the last equality   that
$
\langle {v}, \log^{\kappa}_{p}{P}^{\kappa}_{C}(\gamma(\alpha)) \rangle \leq -({1}/{\alpha})d_{\kappa}^2 (p, {P}^{\kappa}_{C}(\gamma(\alpha))).
$
Since  $\gamma(\alpha):=\mbox{exp}_{p}(-\alpha {v})$, the last inequality is equivalent to \eqref{eq;snlpf}, and this concludes the proof.
\end{proof}
\subsection{Computing the intrinsic $\kappa$-projection into $\kappa$-hyperbolically convex sets}
In this section, we will demonstrate how to explicitly compute the intrinsic $\kappa$-projection for certain $\kappa$-hyperbolically special convex sets. Specifically, for these special $\kappa$-hyperbolically convex sets, we can simplify the computation of the intrinsic $\kappa$-projection on ${C}\subseteq {{\mathbb H}^n_{\kappa}}$ by reducing it to the calculation of the orthogonal Euclidean projection on its associated cone ${K_{C}}$.
\begin{proposition} \label{pr:pcfsT1}
 Let ${C}\subseteq {{\mathbb H}^n_{\kappa}}$ be a  closed $\kappa$-hyperbolically
 convex set, $p:=(p_1, \ldots, p_n, p_{n+1})^T\in{{\mathbb H}^n_{\kappa}}$ and
 $u:=(u_1, \ldots, u_n, u_{n+1})^T\in~{K}_{C}\setminus\{0\}$. Consider the following
 statements:
 \begin{enumerate}
 \item[(i)]
  \begin{equation}\label{eq:SufCondProj1}
	\left(p_{n+1}-u_{n+1}\right)\Big(q^{n+1}-\f{u^{n+1}}{\sqrt{-{\kappa}\lng u, u \rng}}\Big)\ge 0,  \qquad \forall q:=(q_1, \ldots, q_n, q_{n+1})^T\in {C};
\end{equation}
 \item[(ii)]
 \begin{equation} \label{q:SufCondProj2}
 \left(2p_{n+1}e^{n+1}-u\right)^{\top} \Big(q-\frac{u}{\sqrt{-{\kappa}\lng u, u \rng}}\Big)\geq 0, \qquad \forall q\in {C};
 \end{equation}
  \item[(iii)]  If  $ u\in \Pi_{{K}_{C}}(p)$, then
  \begin{equation} \label{q:SufCondProj3}
  {P}^{\kappa}_{C}(p)=\frac{u}{\sqrt{-{\kappa}\lng u, u \rng}}.
 \end{equation}
 Then,  the following implications  hold: (i)$\implies$(ii)$\implies$(iii).
 \end{enumerate}
\end{proposition}
\begin{proof}
To begin, for the sake of simplifying notation, let us consider a given  $u\in{{\mathbb H}^n_{\kappa}}$. We will set
$$
 \eta:=\sqrt{-{\kappa}\langle u, u \rangle}.
 $$
Proof of (i)$\implies$(ii): Let $u$ and $q$ be vectors in $\Hy$. Since $u/\eta$ and $q$ are  in $\Hy$, we can conclude that
\begin{equation} \label{eq:fiscproj}
\Big\langle \frac{u}{\eta}, \frac{u}{\eta}\Big\rangle =-\frac{1}{\kappa}, \qquad \quad \Big\langle \frac{u}{\eta},q\Big\rangle\leq -\frac{1}{\kappa}.
\end{equation}
By combining both the equality and inequality stated in \eqref{eq:fiscproj} and taking into account \eqref{eq:ip}, we derive
\begin{equation} \label{eq:fiscproj1}
(-\J u)\tp \Big(q-\f u{\eta}\Big)= -\eta \Big\langle\frac{u}{\eta},q-\frac{u}{\eta}\Big\rangle\geq 0.
\end{equation}
Since $-u+\J u=2u_{n+1}e^{n+1}$, inequality \eqref{eq:SufCondProj1} can be equivalently stated as
\begin{equation}\label{eq:fiscproj2}
\f12(2p_{n+1}e^{n+1}-u+\J u)\tp\Big(q-\f u{\eta}\Big)\ge 0,
\end{equation}
where $e^{n+1}=(0,\ldots, 0,1)\in {\mathbb R}^{n+1}$. Therefore, by combining \eqref{eq:fiscproj1} with \eqref{eq:fiscproj2}, we obtain  \eqref{q:SufCondProj2}.

\noindent
Proof of (ii)$\implies$(iii): Since ${C}\subset{{\mathbb H}^n_{\kappa}}$ is
a closed $\kappa$-hyperbolically convex set, Proposition~\ref{pr:ccs} implies
that ${K}_{C}$ is a convex cone. Take $q\in {C}$ and $u\in {\mathbb
R}^{n+1}$ such that $u\in \Pi_{{K}_{C}}(p)$. Due to ${\eta}q \in {K}_{C}$ and \eqref{eq:zara}, we have $0\geq (p-u)^{\top} ({\eta}q-u)$. Thus, we can conclude that
\begin{equation}\label{eq:cdepT1}
\Big\langle p, \frac{u}{\eta} \Big\rangle\geq \left\langle p, q\right\rangle- \left\langle p, q\right\rangle+ \Big\langle p, \frac{u}{\eta} \Big\rangle+ \frac{1}{\eta}(p-u)^{\top} ({\eta}q-u).
\end{equation}
On the other hand, by utilizing \eqref{q:SufCondProj2}, we can deduce, after some calculations, that
$$\left\langle p, q\right\rangle+\frac{1}{\eta}\left\langle p, u \right\rangle+ \frac{1}{\eta}(p-u)^{\top} ({\eta}q-u)=\left(p-\J p-u\right)^{\top} \big(q-\frac{1}{\eta}u\big)
=\left(2p_{n+1}e^{n+1}-u\right)^{\top} \big(q-\frac{1}{\eta}u\big)\geq 0, $$
where $e^{n+1}=(0,\ldots, 0,1)\in {\mathbb R}^{n+1}$. Hence, by combining the  last inequality with \eqref{eq:cdepT1}, we have
$
\langle p, \frac{u}{\eta} \rangle\geq \left\langle p, q\right\rangle,
$
for all $q\in {C}$. Thus, by using \eqref{d:prj2} and that $ \eta=\sqrt{-{\kappa}\langle u, u \rangle}$,  we obtain \eqref{q:SufCondProj3}.
\end{proof}
In the following examples  we use  Proposition~\ref{pr:pcfsT1} to present formulas to intrinsic $\kappa$-projection for certain special $\kappa$-hyperbolically convex sets.
\begin{example}\label{ex:simp}
Let ${K}=\cone\left\{v^1,\dots,v^n,v^{n+1}\right\}$ be a simplicial cone, where
$(v^i)\tp e^{n+1}=0$ with  $v^i\neq 0$, for all  $i=1,\dots,n$  and
$v^{n+1}=e^{n+1}$. Consider a $\kappa$-hyperbolically convex set  ${C}$ and  corresponding cone spanned by  the $\kappa$-hyperbolically convex set ${C}$  given, respectively,  by
\begin{equation*}
	{C}:= {K} \cap {{\mathbb H}^n_{\kappa}},  \qquad \qquad {K_{C}}:=\left\{ p:~p\in {\rm Int}({L})\cap {K}
	\right\}\cup\{0\}.
\end{equation*}
We notice that ${K_{C}}$   is  a convex cone. If  $u\in\Pi_{{K}_{C}}(p)$, then the intrinsic $\kappa$-projection into the $\kappa$-hyperbolically convex set ${C}$ is given by
 \begin{equation}  \label{eq:ProjCircSet}
  {P}^{\kappa}_{C}(p)=\frac{u}{\sqrt{-{\kappa}\lng u, u \rng}}.
 \end{equation}
Indeed, take   $p:=(p_1, \ldots, p_n, p_{n+1})\in{{\mathbb H}^n_{\kappa}}$ and $u:=(u_1, \ldots, u_n,
 u_{n+1})\in~{C}$.   If  $u\in\Pi_{{K}_{C}}(p)$, then  $u_{n+1}=p_{n+1}$. In fact, suppose, on the contrary, that $u_{n+1}\ne p_{n+1}$.   Since
 $u\in\Pi_{{K}_{C}}(p)$ and $u_{n+1}\ne p_{n+1}$, we infer that  $p\notin {K}_{C}$ and $u$ is on the boundary of ${K}$. Hence, there exists $j\in\{1,\dots,n,n+1\}$, such  that
 \begin{equation} \label{eq:face}
 u\in \mathcal{F}:=\cone\left\{\left\{v^1,\dots,v^n,v^{n+1}\right\}\setminus\left\{v^j\right\}\right\}.
 \end{equation}
 As  $u\in\Pi_{{K}_{C}}(p)$, we infer that $u$ belongs to the boundary of ${K}$. Hence,  $u$ is on one of the $(n-1)$-dimensional faces of ${K}$.  If $j=n+1$, then
$u\in\cone\left\{v^1,\dots,v^n\right\}$. Thus,  $u\tp e^{n+1}=0$. Since, $e^{n+1}\in{K}_{C}$, it follows from \eqref{eq:zara} that $p_{n+1}=(p-u)\tp e^{n+1}\le 0$, but this is
absurd, because
$p\in\Hy$. Hence, it follows that $j\in\{1,\dots,n\}$. By the definition of the projection and \eqref{eq:face}, it follows
that $p-u$ must be perpendicular to the face $\mathcal{F}$.  On the other hand,  due  to $j\in\{1,\dots,n\}$,
it follows that $e^{n+1}=v^{n+1}\in\mathcal{F}$. Hence,
$p_{n+1}-u_{n+1}=(p-u)\tp e^{n+1}=0$ and we have $u_{n+1}=p_{n+1}$. As
consequence, item (i) of Proposition~\ref{pr:pcfsT1} holds, i.e.,
\eqref{eq:SufCondProj1} holds. Therefore, by using item (iii) of Proposition~\ref{pr:pcfsT1},  we obtain  that the projection of $p$ into ${C}$ is in fact  given by \eqref{eq:ProjCircSet}.
 \end{example}
 In the next example, we utilize Example~\ref{ex:simp} to present an explicit formula for computing the intrinsic $\kappa$-projection into the following $\kappa$-hyperbolically convex set
 \begin{equation} \label{eq:hpc}
 {C}_{+}=\{p\in {\mathbb H}^n_{\kappa}:~ p\geq 0\},
 \end{equation}
 involving only nonnegative constraints, which are fundamental in many applications, see \cite{Chen2010}.
 \begin{example}\label{ex:PositOrthat}
  Let us compute a formula for $\kappa$-projecting intrinsically  a point into  ${C}_{+}$ given by \eqref{eq:hpc}.
  To achieve that, it is important to note that the cone spanned by the set ${C}_{+}$ can be expressed as
\begin{equation} \label{eq:pccone2}
	{{K}_{{C}_{+}}}:=({\mathbb R}^{n+1}_{+}\cap{\rm Int}({L}))\cup\{0\},
\end{equation}
 which is a convex cone.  It follows from  ${{K}_{{C}_{+}}} \subseteq  {\mathbb R}^{n+1}_{+}$ and \eqref{eq:pirn} that
\begin{equation} \label{eq:appsb}
{\inf}_{{z\in ({\rm Int}({L})\cap {\mathbb R}^{n+1}_{+})\cup\{0\}}}\|p-z\|\geq {\inf}_{{z\in {\mathbb R}^{n+1}_{+}}}\|p-z\|.
\end{equation}
On the other hand, we know that
\begin{equation} \label{eq:appsb2}
\{p^{+}\}=\Pi_{{\mathbb R}^{n+1}_{+}}(p)={\arg\min}_{{z\in {\mathbb R}^{n+1}_{+}}}\|p-z\|.
\end{equation}
The combination of  \eqref{eq:cc}  for $\alpha=1$ with  \eqref{eq:hs} implies  that $p=(p_1,
\ldots, p_n,p_{n+1})\in {\rm Int}({L})$. Hence, considering the fact that
$p_{n+1}\geq \frac{1}{\sqrt{\kappa}}>0$, we obtain that
$$
p^{+}_{n+1}=p_{n+1}> \sqrt{p_1^2+ \cdots +p_n^2}\geq  \sqrt{(p_1^+)^2+ \cdots +(p_n^+)^2},
 $$
 where $p^{+}:=(p^{+}_1, \ldots, p^{+}_n,p^{+}_{n+1})$ and $p^{+}_i:=\max\{p_i,
 0\}$ for $i=1, \ldots, n, n+1$. Thus, we have  $p^{+}\in  ({\rm Int}({L})\cap {\mathbb R}^n_{+})\cup\{0\} $, which combined with  \eqref{eq:appsb} and \eqref{eq:appsb2}  yields  $p^{+}=\Pi_{{K}_{{C}_{+}}}(p).$ Hence, as  ${\mathbb R}^n_{+}=\cone\left\{e^1,\dots,e^n,e^{n+1}\right\}$ and ${{C}_{+}}:={\mathbb R}^n_{+} \cap {{\mathbb H}^n_{\kappa}}$,  it follows from   Example~\ref{ex:simp} that   the projection of $p\in {{\mathbb H}^n_{\kappa}}$ into the set  ${C}$ is given by
\begin{equation*}
  {P}^{\kappa}_{{C}_{+}}(p)=\frac{p^{+}}{\sqrt{-{\kappa}\lng p^{+}, p^{+} \rng}}.
\end{equation*}
 \end{example}

 \begin{example}\label{ex:cc}
Let ${C}:= {L}_{\alpha} \cap {{\mathbb H}^n_{\kappa}}$, where ${L}_{\alpha}$ is  a circular cone and $\alpha\ge 1$ is a real constant. Since
$\alpha\ge 1$, we have  ${K_{C}}=\Lo_\alpha$  and it is a closed convex cone.  If $p\in\Hy$ and $u=\Pi_{{K}_{C}}(p)$, then $u\ne 0$ and \eqref{eq:SufCondProj1} holds.  Indeed, since $e^{n+1}\in{K}_{C}$ and $u=\Pi_{{K}_{C}}(p)$, it follows from \eqref{eq:zara} that
		\begin{equation}\label{eq:n+1}
			p_{n+1}-u_{n+1}=(p-u)\tp e^{n+1}\le 0.
		\end{equation}
		On the other hand, since $q\in {C}={K}_{C}\cap\Hy$, we have
		$q_{n+1}^2\ge\alpha^2(q_1^2+\dots+q_n^2)=\alpha^2 (q_{n+1}^2-\f1\kappa),$
		which after some algebraic manipulations implies  that
		\begin{equation}\label{eq:qn+1}
			q_{n+1}\le\f\alpha{\sqrt{\kappa(\alpha^2-1)}},
		\end{equation}
		with equality holding if and only if $q$ is on the boundary of $\C$.
		If $p\in\C$, then the projection formula is trivial. Suppose that
		$p\notin\C$. Thus, $u$  is on the boundary of $\K_{C}$, which implies that
		the point
		$$
		\frac{u}{\sqrt{-{\kappa}\langle u, u \rangle}}
		$$
		 is on the boundary of the set $\C$. By  using the same
		 argument as the one used to show \eqref{eq:qn+1}, we  can  show  that
		$
			{u_{n+1}}/{\sqrt{-{\kappa}\langle u, u \rangle}}=\alpha/{\sqrt{\kappa(\alpha^2-1)}}.
		$
		Thus,  the last equality together with   \eqref{eq:qn+1} implies  that
		\begin{equation}\label{eq:diffsign}
			q_{n+1}-\frac{u_{n+1}}{\sqrt{-{\kappa}\langle u, u \rangle}}\leq 0.
		\end{equation}
		It follows from \eqref{eq:n+1} and \eqref{eq:diffsign} that \eqref{eq:SufCondProj1} holds. Therefore,  item (iii) of Proposition~\ref{pr:pcfsT1} implies that
	\begin{equation}\label{eq:eunelo}
  		{P}^{\kappa}_{C}(p)=\frac{u}{\sqrt{-{\kappa}\lng u, u \rng}}.
	\end{equation}
\end{example}
\subsection{The Lorentz projection into convex sets}

In this section, we will introduce the Lorentz projection and explore its connections with the Euclidean orthogonal projection and the intrinsic $\kappa$-projection.
\begin{definition}
	Let $\K\subseteq\R^{n+1}$ be a convex set, $x\in\R^{n+1}$ and $y\in\K$. The point
$y$ is called a \emph{Lorentz projection} (or shortly \emph{L-projection}) of $x$ into $\K$ if $\lng x-y,z-y\rng\le 0$, for all
$z\in\K$.
\end{definition}
Denote by $\Lambda_\K(x)$ the set of {\it L-projections} of $x\in\R^{n+1}$ into $\K$, i.e.,
	\begin{equation}\label{eq:L-Proj}
	\Lambda_\K(x):=\{ y\in\K:~ \lng x-y,z-y\rng\le 0, ~\forall z\in\K\}.
	\end{equation}
As the Lorentzian inner product is not positive definite, given $x\in\R^{n+1}$,
the set $\Lambda_\K(x)$ might be empty or even unbounded.  However, the general
analysis of $\Lambda_\K(x)$ is out of the scope of this paper. As we will see
later, if $p\in\Hy$ and $\K\subseteq\Lo$ is a cone such that $\K\cap\Hy$ is a
closed  $\kappa$-hyperbolically convex set, then $\Lambda_\K(p)$ is nonempty and
contains precisely one non-zero vector in addition to the zero vector. In the
next remark we show that zero vector belongs to the L-projections.
\begin{remark}  \label{re:0inlp}
If $\K\subseteq\Lo$ such that $\K\cap\Hy$ is a $\kappa$-hyperbolically
convex set and $p\in\Hy$, then   it follows from $p\in\Hy\subseteq \inte(\Lo)$,  \eqref{eq:L-Proj} and Remark~\ref{re;PLC} that $0\in 	\Lambda_\K(p)$, for all $p\in {\mathbb H}^{n}_{\kappa}$.
\end{remark}
In the following lemma, we present equivalent conditions for a point to be an L-projection of a point into a cone.  They are corresponding to relations \eqref{eq:zara} fulfilled by the Euclidean orthogonal projection.
\begin{lemma}\label{lem:lmp}
	Let $\K\subseteq\R^{n+1}$ be a convex cone, $x\in\R^{n+1}$ and $y\in\K$.
	Then, $y\in\Lambda_\K(x)$ if and only if
	the following two relations hold:
	\begin{equation}\label{eq:lmp}
		\lng x-y,z\rng\leq 0,\mbox{ }  \forall\, z\in \K\qquad\mbox{and} \qquad \lng x-y,y\rng=0.
	\end{equation}
\end{lemma}
\begin{proof}
	We have that $y$ is an L-projection of $x$ into $\K$ if and only if $\lng x-y,z-y\rng\le 0$,
	for all $z\in\K$. Hence, we must prove that $\lng x-y,z-y\rng\le 0$, for all $z\in\K$,  is equivalent
	to \eqref{eq:lmp}. Suppose that $\lng x-y,z-y\rng\le 0$ holds. Then,
	\(
		\lng x-y,\lambda z-y\rng\le 0,
	\)
	for all $z\in\K$ and $\lambda>0$. By letting $\lambda=2$ and $z=y$, we obtain
	$\lng x-y,y\rng\le 0$, while by letting $\lambda\to 0$ we obtain
	$\lng x-y,y\rng\ge 0$. Hence, $\lng x-y,y\rng=0$, which
	inserted into $\lng x-y,z-y\rng\le 0$ also implies $\lng x-y,z\rng\le 0$. Therefore, $\lng x-y,z-y\rng\le 0$ holds, for all $z\in\K$.  Next, suppose
	that \eqref{eq:lmp} holds. Then, $\lng x-y,z-y\rng\le 0$, for all $z\in\K$,  follows from the
	bilinearity of the Lorentzian   inner product.
\end{proof}
In the  following  corollary we use Lemma~\ref{lem:lmp} to compute a  L-projection into a hyperplane.
\begin{corollary}\label{cor:lph}
Let $a\in\R^n$   be such that $\lng a,a\rng\neq 0$ and  $\V:=\{x\in\R^n:~a\tp x=0\}$ be a hyperplane. Then
$$
\Lambda_{\V}(y)\setminus\{0\}\ni  y-\f{\lng y,\J a\rng}{\lng a,a\rng}\J a, \qquad \forall  y\in\R^{n+1}.
$$
\end{corollary}
\begin{proof}
Since $\V$  is a hyperplane,  its follows from Lemma~\ref{lem:lmp}   that  $u\in \Lambda_{\V}(y)\setminus\{0\}$  if and only if the following equality holds
\begin{equation}\label{eq:pslf2}
(\J y-\J u)\tp x=\lng y-u,x\rng= 0, \qquad \forall x\in \V.
 \end{equation}
 Since $a$ is the normal of the hyperplane $\V$, it follows from \eqref{eq:pslf2} that  there exists $t\in \R$ such that $\J y-\J u=ta$, or equivalently,
\begin{equation}\label{eq:pslf3}
 u=y- t\J a
 \end{equation}
 Using \eqref{eq:pslf3} we obtain that $\langle u,\J a\rangle= \langle y,\J
 a\rangle - t \langle \J a,\J a\rangle$. Since  $u\in \V$, we have $\langle u,\J
 a\rangle=0$, and taking into account that   $\langle \J a, \J a\rangle=\langle a, a\rangle \neq 0$, we conclue that
$t={\langle y,\J a\rangle}/{\langle  a, a\rangle}.$
Substituting $t$ given by the last equality into \eqref{eq:pslf3} we obtain the desired result.
\end{proof}
In the following, we will revisit the concept of a linear cone-complementarity
problem and explore its potential connection to \eqref{eq:lmp} in a particular
case. To establish this relationship, let us define
$
\K^*:=\{w\in \R^{n+1}:~w\tp v\ge 0,~\forall v\in \K\},
$
the {\it dual cone} of the cone $\K\subseteq\R^{n+1}$.
\begin{definition}
	Let $\K\subseteq\R^{n+1}$ be a closed convex cone, $M\in\R^{(n+1)\times (n+1)}$ be a
	matrix and $q\in\R^{n+1}$ be a vector. Then, the \emph{linear cone-complementarity
	problem} LCP$(q,M,\K)$ defined by $q$, $M$ and $\K$ is the problem of
	finding an $y\in\K$ such that $z:=q+My\in\K^*$ and $y\tp z=0$.
\end{definition}
In the upcoming corollary, we establish a connection between the L-projection and
a particular linear cone-complementarity problem.
To achieve that, denote by $\ovl{S}$ the {\it topological closure} of a set ${S}\subseteq\R^{n+1}$.
\begin{corollary}\label{cor:lcp}
	Let $\K\subseteq\R^{n+1}$ be a convex cone. Then, $y\in\Lambda_\K(x)$
	if and only if $y \in \K$ and it is a solution of the linear cone-complementarity
	problem LCP$(-\J x,\J,\ovl\K)$.
\end{corollary}
\begin{proof}
	From Lemma \ref{lem:lmp} we obtain that  $y\in\Lambda_\K(x)$ if and only
	if $\lng x-y,y\rng=0$ and  $\lng x-y,z\rng\le 0$, for any $z\in\K$. The
	latter two formulas are equivalent to $(\J y-\J x)\tp y=0$
	and $(\J y-\J x)\tp z\ge 0$, for any
	$z\in\K$, respectively. Hence, $y\in\Lambda_\K(x)$ if and only if $(\J y-\J x)\tp y=0$ and
	$\J y-\J x\in\K^*$. Therefore, $y\in\Lambda_\K(x)$
	if and only if it is a solution of the linear cone-complementarity
	problem LCP$(-\J(x),\J,\K)$.
\end{proof}
In the subsequent theorem, we demonstrate that every nonzero L-projection gives rise to an intrinsic $\kappa$-projection, which can be obtained by solving the linear complementarity problem presented in the previous corollary.
\begin{theorem}\label{th:hypp}
Let ${C}\subseteq {{\mathbb H}^n_{\kappa}}$ be a  closed  $\kappa$-hyperbolically convex set, $p\in{{\mathbb H}^n_{\kappa}}$ and $u\in\Lambda_{\K_\C}(p)\setminus\{0\}$. Then,
	\begin{equation}\label{eq:hypp}
  		{P}^{\kappa}_{C}(p)=\frac{u}{\sqrt{-{\kappa}\lng u, u \rng}}.
 	\end{equation}

\end{theorem}
\begin{proof}
Let $p\in{{\mathbb H}^n_{\kappa}}$ and $u\in\Lambda_{\K_\C}(p)\setminus\{0\}\subseteq {L}$. Considering that   $u\in\Lambda_{\K_\C}(p)\setminus\{0\}\subseteq \inte{L}$,   we have ${u}/{\sqrt{-\kappa \lng u,u\rng}}\in\Hy.$ Hence,   $\langle u/\sqrt{-\kappa \langle u,u\rangle}, q\rangle
\leq -1/\kappa$, for any $q\in\Hy$. Thus, we conclude that
\begin{equation}\label{eq:cipicl}
  \lng u, q\rng \leq -\frac{1}{\kappa}\sqrt{-{\kappa}\lng u, u \rng}<0,\qquad \forall q\in \Hy.
\end{equation}
On the other hand, the  conditions  \eqref{eq:lmp}  in Lemma~\ref{lem:lmp} imply, respectively,  that
\begin{equation}\label{eq:cipf}
\lng p, q\rng\le \lng u,q\rng, \qquad  \lng p, u\rng= \lng u, u\rng, \qquad \forall q\in\C\subseteq\K_\C.
\end{equation}
By using the  inequality in \eqref{eq:cipf}  together with the second inequality
in \eqref{eq:cipicl}, we conclude that
$
		\lng p,q\rng\le\lng u,q\rng  <0,
$
for all $q\in\C$. Since the second inequality in \eqref{eq:cipicl}  yields $ \lng u, q\rng^2
\geq -(1/\kappa)\lng u, u \rng$, we obtain from the  last inequality and the equality in \eqref{eq:cipf} that
\begin{equation*}
		\lng p,q\rng^2\ge\lng u,q\rng^2\ge-\f1\kappa\lng u,u\rng=-\f1\kappa\lng p,u\rng=\f{\lng
		p,u\rng^2}{-\kappa\lng u,u\rng}=\Big\langle p, \frac{u}{\sqrt{-\kappa \langle u,u\rangle}}\Big\rangle^2, \qquad \forall q\in\C.
\end{equation*}
Since  $p,q \in\Hy$ implies  $\lng p,q\rng<0$, it follows  from the last inequality  that
\begin{equation*}
		\langle p,q\rangle \leq  \big\langle p, \frac{u}{\sqrt{-\kappa \langle u,u\rangle}}\big\rangle, \qquad \forall q\in\C.
\end{equation*}
Therefore, by using \eqref{d:prj2}, we deduce the desired equality in \eqref{eq:hypp}
.
\end{proof}
The following theorem establishes an equivalence between the intrinsic $\kappa$-projection and the L-projection. In particular, it shows  that the L-projection contains a nonzero vector.
\begin{theorem}\label{th:lp}
	Let ${C}\subseteq {{\mathbb H}^n_{\kappa}}$ is a closed  $\kappa$-hyperbolically convex set, $p\in{{\mathbb H}^n_{\kappa}}$ and $v\in\C$. Then, the following
	statements hold:
	\begin{enumerate}[(i)]
		\item{\label{pc}}\hspace{0mm} We have $v={P}^{\kappa}_{C}(p)$, if and only if
		$-\kappa\lng p,v\rng v\in\Lambda_{\K_\C}(p)\setminus\{0\}.$
	\item\label{ex}\hspace{0mm} The set $\Lambda_{\K_\C}(p)\setminus\{0\}\ne\varnothing$.
	\end{enumerate}
\end{theorem}
\begin{proof}
{\it Item~(i)}:  Suppose $v={P}^{\kappa}_{C}(p)$ and let $u:=-\kappa\lng p,v\rng v$. Since $p,v\in\Hy$, it follows that $-\kappa\lng p,v\rng>0$ and $\lng p+\kappa\lng p,v\rng v,v\rng=0$. Hence, we conclude that
	\begin{equation} \label{eq:vfp}
	 u\in\K_\C\setminus\{0\}, \qquad \lng p-u,u\rng=0.
	\end{equation}
	On the other hand, for any $z\in\K_\C$, we have
	\begin{equation} \label{eq:vfpn1}
	q:=\f z{\sqrt{-\kappa\lng z,z\rng}}\in\C.
	\end{equation}
Thus, by combining \eqref{eq:Intdist} with Lemma~\ref{le:CosLaw1} and denoting
$\theta_v$ as the angle between the vectors $\log^{\kappa}_{{v}}{p}$ and
$\log^{\kappa}_{v}{q}$, we can deduce, after performing some calculations, that
	\begin{equation} \label{eq:vfpn}
	\kappa \lng p,q\rng+\kappa^2\lng v, p\rng\lng v,q\rng=\sqrt{1- \kappa^2\lng v, p\rng^2}\sqrt{1- \kappa^2\lng v, q\rng^2}\cos\theta_v.
	\end{equation}
Since $v={P}^{\kappa}_{C}(p)$ and $\theta_v$ is the angle between the vectors
$\log^{\kappa}_{{v}}{p}$ and $\log^{\kappa}_{v}{q}$, it follows from
\eqref{d:uexp} that 	$\cos\theta_v\leq 0$. Thus, by using \eqref{eq:vfpn}, we
infer that
$
\kappa \lng p,q\rng+\kappa^2\lng v, p\rng\lng v,q\rng\leq 0.
$
Hence,  \eqref{eq:vfpn1},  $u=-\kappa\lng p,v\rng v$ and the  last inequality imply that
$$
\lng p-u,z\rng =\frac{1}{\kappa}\sqrt{-\kappa\lng z,z\rng}\lf(\kappa \lng p,q\rng+\kappa^2\lng p,v\rng\lng v,q\rng\rg)\leq 0, \qquad \forall z\in\K_\C.
$$
Therefore, the latter  inequality, \eqref{eq:vfp} and  Lemma \ref{lem:lmp} imply
that $u\in\Lambda_{\K_\C}(p)\setminus\{0\}$, and the proof is concluded.

Conversely, suppose that $u:=-\kappa\lng p,v\rng
v\in\Lambda_{\K_\C}(p)\setminus\{0\}$. Since
$v\in\C\subseteq\Hy$, it follows that $\lng v,v\rng=-1/\kappa$. Hence, the
results follow from Theorem \ref{th:hypp} after performing some calculations.

\noindent
{\it Item~(ii)}:  Since ${C}\subseteq {{\mathbb H}^n_{\kappa}}$ is closed $\kappa$-hyperbolically
	 convex set, for each  $p\in{{\mathbb H}^n_{\kappa}}$ there exists $v={P}^{\kappa}_{C}(p)$. Therefore, the result is a direct consequence of item~\ref{pc}.
\end{proof}

In the following corollary, we assume that ${K}\subseteq{L}$ is a convex cone
such that ${K}\cap\Hy$ is a closed $\kappa$-hyperbolically convex set, and $0$
belongs to ${K}$. We will prove various properties of the L-projection of a point
$p\in\text{int}({L})$ into ${K}$.

\begin{corollary}\label{cor:bd}
	Let $p\in\inte(\Lo)$,  $\K\subseteq\Lo$ be a convex cone,
	such that $\C:=\K\cap\Hy$ is a closed $\kappa$-hyperbolically convex set and $0\in\K$.
	Then, the following statements hold:
	\begin{enumerate}[(i)]
		\item\label{unique}\hspace{0mm} $0\in \Lambda_\K(p)$ and
			$\Lambda_\K(p)\setminus\{0\}$ has
			exactly one element.
		\item\label{bd}\hspace{0mm} If $p\notin\K$, then $\Lambda_\K(p)\subseteq\p\K$, where $\p$ denotes boundary.
		\item\hspace{0mm} If $p\in\K$, then $\Lambda_\K(p)=\{0,p\}$.
		\item\label{natural}\hspace{0mm} $ u=(u_1, \dots, u_n, u_{n+1})^T\in\Lambda_\K(p)\iff
			u=\Pi_{\bar\K}\lf(\J p+2u_{n+1}e^{n+1}\rg)\iff
			u=\Pi_{\tau_u(\bar\K)}(p),$ where  \[\tau_u(\bar{K})
			:=\left\{\tau_u(x):~x=(x_1, \dots, x_n, x_{n+1})^T\in\bar{K}
		\right\}\]
			and $\tau_y:\R^{n+1}\to\R^{n+1}$ is defined by  $\tau_y(z):=z+2(y_{n+1}-z_{n+1})e^{n+1}$ with  $y=(y_1, \dots, y_n, y_{n+1})^T$ and $z=(z_1, \dots, z_n, z_{n+1})^T.$
		\item\label{eu=lo}\hspace{0mm} We have
			$\varnothing\ne\Lambda_{\K}(p)\setminus\{0\}=\Pi_{\K}(p)$, if and
			only if
			\[\lf(\Lambda_{\K}(p)\setminus\{0\}\rg)
			\cap\{x\in\R^{n+1}:~x_{n+1}=p_{n+1}\}\ne\varnothing.\]
		\item\label{solve}\hspace{0mm} Finding an $u=(u_1, \dots, u_n, u_{n+1})^T\in\Lambda_\K(p)$ is equivalent to determining
			$u_{n+1}$ from any of the  following  equations
			\[u_{n+1}=\lf(e^{n+1}\rg)\tp\Pi_{\bar\K}\lf(\J
			p+2u_{n+1}e^{n+1}\rg),\qquad \qquad
			u_{n+1}=\lf(e^{n+1}\rg)\tp\Pi_{\bar\K+2\lf[u_{n+1}-\lf(e^{n+1}\rg)\tp
			\bar\K\rg]e^{n+1}}(p),\]
			followed by using any of the following equations
			\[u=\Pi_{\bar\K}\lf(\J p+2u_{n+1}e^{n+1}\rg),\qquad \qquad
			u^{n+1}=\Pi_{\bar{K}+2\lf[u_{n+1}-\lf(e^{n+1}\rg)\tp\bar{K}\rg]e^{n+1}}(p),\]
			where ${\bar{K}+2\lf[u_{n+1}-(e^{n+1})\tp\bar{K}\rg]e^{n+1}}
			:=\left\{x+2(u_{n+1}-(e^{n+1})\tp x)e^{n+1}:~ x\in\bar{K}\right\}$.
		\item\label{usolve}\hspace{0mm} Finding the $u=(u_1, \dots, u_n, u_{n+1})^T\in\Lambda_\K(p)\setminus\{0\}$ is equivalent to determining $u_{n+1}\ne 0$ from the  following equation
			\[u_{n+1}=\lf(e^{n+1}\rg)\tp\Pi_{\bar\K}\lf(\J
				p+2u_{n+1}e^{n+1}\rg),\]
				followed by using the equation
				$u=\Pi_{\bar\K}\lf(\J p+2u_{n+1}e^{n+1}\rg).$
		\item\hspace{0mm} For any $\lambda>0$ we have $\Lambda_\K(\lambda
			p)=\lambda\Lambda_\K(p).$
		\item\label{pp}\hspace{0mm} If
			$\Lambda_{\K}(p)\setminus\{0\}=\lambda\Pi_{\K}(p)$,
			for some $\lambda>0$, then
			\begin{eqnarray*}
				\Lambda_K(p)\setminus\{0\}=
				\f{\lng p,\Pi_\K(p)\rng}{\lng\Pi_\K(p),\Pi_\K(p)\rng}
				\Pi_\K(p).
			\end{eqnarray*}
		\item\hspace{0mm} Suppose that $\Pi_\K(p)\ne\varnothing$. We have
			$\Lambda_{\K}(p)\setminus\{0\}=\lambda\Pi_{\K}(p)$,
			for some $\lambda>0$ if and only if
			\begin{equation} \label{eq:rslpep}
				\Pi_\K(p)
				=\argmin\Big\{\f{\lng p,z\rng}{\lng
				\Pi_\K(p),z\rng}:z\in\K\setminus\{0\}\Big\}.
			\end{equation}
		\item\hspace{0mm} Let $Q\in\R^{(n+1)\times (n+1)}$ with
				$Q\tp\J Q=\J$. Then, the set $Q\C$ is $\kappa$-hyperbolically
				convex and
			\[\Lambda_{\K_{Q\C}}(p)=\Lambda_{Q\K}(p)=\Lambda_\K\big(Q\tp\J
			p\big).\]
	\end{enumerate}
\end{corollary}

\begin{proof}
{\it Item~(i)}:  We will use the following notations: For  $\kappa >0$ and $a\in\Lo\setminus\{0\}$ we denote
\[
\eta(a):=\sqrt{-\kappa\lng a,a\rng}, \qquad \quad  \mu(a):=\f a{\eta(a)}\in\Hy.
\]  It follows from Remark~\ref{re:0inlp} that $0\in\Lambda_\K(p)$. Let
$p\in\inte(\Lo)$ and $ {{\mathbb H}^n_{\kappa}}$ be the $\kappa$-hyperbolic space
form  such that  $\kappa:=-1/{\lng p,p\rng}.$ We have $\K=\K_\C$. Since ${C}$ is a closed $\kappa$-hyperbolically convex set,  ${P}^{\kappa}_{C}(p)$ contains
exactly one element. Let  ${\hat u},{\tilde u}\in\Lambda_\K(p)\setminus\{0\}$.
Then, Theorem \ref{th:hypp} implies that $\mu({\hat u})=\mu({\tilde
u})={P}^{\kappa}_{C}(p)$. Hence, ${\tilde u}=\lambda {\hat u}$, where
$\lambda:={\eta({\tilde u})}/{\eta({\hat u})}>0.$ On the other hand,
Lemma~\ref{lem:lmp} yields $\lng {\hat u}-p,{\hat u}\rng=0$ and  $\lng p-\lambda
{\hat u},{\hat u}\rng=(1/\lambda)\lng p-{\tilde u},{\tilde u}\rng=0.$ By summing
up the last two inequalities, we get $(1-\lambda)\lng {\hat u},{\hat u}\rng=0$.
Therefore, based on Remark~\ref{re;PLC}, which implies that $\langle {\hat u},{\hat u}\rangle <0$, we conclude that $\lambda=1$.
Hence, ${\hat u}={\tilde u}$. Hence, $\Lambda_\K(p)\setminus\{0\}$ contains at most one element. Therefore, item \ref{ex} of Theorem \ref{th:lp} implies that $\Lambda_\K(p)\setminus\{0\}$ contains exactly one element.

\noindent
 {\it Item~(ii)}:  Let   $p\in\inte(\Lo)$ and $ {{\mathbb H}^n_{\kappa}}$ be the $\kappa$-hyperbolic space form  such that  $\kappa:=-1/{\lng p,p\rng}.$ We have $\K=\K_\C$, where $\C=\K\cap\Hy$. Since $p\notin \K$ it follows that $p\notin\C$. This, together with $\C$ being a closed $\kappa$-hyperbolically convex set, implies ${P}^{\kappa}_{C}(p)\in\p\C$.  Hence, formula \eqref{eq:hypp} of Theorem \ref{th:hypp} implies that $\Lambda_\K(p)\subseteq\p\K$.

\noindent
 {\it Item~(iii)}: 	From the definition of the L-projection and $p\in\Lo$ it follows that $p\in\Lambda_\K(p)$. Hence, the result is an immediate consequence of \ref{unique}.

\noindent
 {\it Item~(iv)}: The first equivalence follows from Corollary \ref{cor:lcp}, because the right hand side of the equivalence is the natural equation of LCP$\lf(-\J p,\J,\ovl\K\rg)$, see \cite[Proposition~1.5.8]{FacchineiPangI}. The second equivalence follows from the formula $\lng y-x,z-x\rng=(y-x)\tp\lf(\tau_y(z)-x\rg)$.

\noindent
 {\it Item~(v)}: Suppose tha $\varnothing\ne\Lambda_\K(p)\setminus\{0\}=\Pi_\K(p)$ and  let $u\in\Lambda_\K(p)\setminus\{0\}=\Pi_\K(p)$. Then, Lemma~\ref{lem:lmp} and equation \eqref{eq:zara} imply  $0=\lng p-u,u\rng-(p-u)\tp u=2\lf(p_{n+1}-u_{n+1}\rg)u_{n+1}=0$.  Hence, $u_{n+1}=0$ or $u_{n+1}=p_{n+1}$. If $u_{n+1}=0$, then item \ref{natural} implies $u=\Pi_\K(\J p)=0$, where the latter equality follows from $\J p\in\J\K\subseteq\J\Lo=-\Lo=-\Lo^*$. Hence, $u=0$, which is a contradiction. Thus, $u_{n+1}\ne 0$ and therefore $u_{n+1}=p_{n+1}$.  Conversely, suppose that $u\in\{x\in\R^{n+1}:x_{n+1}=p_{n+1}\}$. Then, $(p-u)\tp (u-z)=\lng p-u,u-z\rng\le 0$,  for all $z\in\K$.  Therefore, $u\in\Lambda_{\K}(p)\setminus\{0\}$ if and only if $u\in\Pi_{\K}(p)$.

\noindent
 {\it Item~(vi)}: It is an immediate consequence of \ref{natural} and \ref{unique}.

 \noindent
 {\it Item~(vii)}: It is an immediate consequence of \ref{solve} and $\Pi_\K(\J
 p)=0$, which follows from \[\J p\in\J \Lo=-\Lo=-\Lo^*=-\lf(\inte(\Lo)\cup\{0\}\rg)^*\subseteq-\K^*.\]

 \noindent
 {\it Item~(viii)}: It follows from Lemma \ref{lem:lmp}.

 \noindent
 {\it Item~(ix)}: Suppose that
			$\Lambda_{\K}(p)\setminus\{0\}=\lambda\Pi_{\K}(p)$,
			for some $\lambda>0$ and let $v\in\Pi_\K(p)$, which
			exists and it is unique because of our assumption and item
			\ref{unique}. Then,
			$\lambda v\in\Lambda_\K(p)$. Hence, it follows from Lemma
			\ref{lem:lmp} that $\lng p-\lambda v,v\rng=0$, which
			implies $\lambda={\lng p,v\rng}/{\lng v,v,\rng }.$
			Therefore,
			\begin{equation*}
				\Lambda_\K(p)\setminus\{0\}
				=\f{\lng p,\Pi_\K(p)\rng}{\lng\Pi_\K(p),\Pi_\K(p)\rng}
				\Pi_\K(p),
			\end{equation*}

 \noindent
 {\it Item~(x)}:  Suppose that
			$\Lambda_{\K}(p)\setminus\{0\}=\lambda\Pi_{\K}(p)$,
			for some $\lambda>0$ and let $v\in\Pi_\K(p)$, which
			exists and it is unique because of our assumption and item \ref{unique}.
			Then, item \ref{pp}
			implies that $\lambda={\lng p,v\rng}/{\lng v,v\rng
			}.$ Then, by using Lemma
			\ref{lem:lmp}, we have $\lng p-\lambda v,z\rng\le 0$, for
			any $z\in\K\setminus\{0\}$. Thus, by using Remark~\ref{re;PLC},
			it follows that \[\f{\lng p,z\rng}{\lng v,z\rng}\ge\lambda
			=\f{\lng p,v\rng}{\lng v,v\rng}, \qquad \forall z\in\K\setminus\{0\}, \]
			which implies \eqref{eq:rslpep}. Conversely, suppose that \eqref{eq:rslpep} holds: Let $v\in\Pi_\K(p)\ne\varnothing$, which is unique
			because $\K$ is convex. Denote
			\begin{equation}\label{eq:lambda}
				\lambda:=\f{\lng p,v\rng}{\lng v,v\rng}.
			\end{equation}
			Then,
			$\lambda>0$, because of Remark \ref{re;PLC}. Let any
			$z\in\K$. By using equation \eqref{eq:lambda}, we obtain that
			\begin{equation}\label{eq:lpc}
				\lng p-\lambda v,v\rng=0,
			\end{equation}
			Due to equation \eqref{eq:rslpep}, we also have  $({\lng p,z\rng}/{\lng v,z\rng})\ge ({\lng p,v\rng}/{\lng v,v\rng})$, which implies
			\begin{equation}\label{eq:lpi}
				\lng p-\lambda v,z\rng\le 0.
			\end{equation}
			Hence, by using \eqref{eq:lpc}, \eqref{eq:lpi}, Lemma
			\ref{lem:lmp} and item \ref{unique},
			it follows that $\lambda v=\Lambda_\K\setminus\{0\}$.
			In conclusion,
			$\Lambda_{\K}(p)\setminus\{0\}=\lambda\Pi_{\K}(p),$
			with $\lambda>0$.

	\noindent
	{\it Item~(xi)}: Since for any $x,y\in\R^{n+1}$ we have $\lng Qx,Qy\rng=\lng x,y\rng$ and
	$Q\Hy=\Hy$, it follows that $Q$ is an isommetry of $\Hy$ and $Q^{-1}=JQ\tp J$. Hence, $Q\C$ is hyperbolically convex
	and $\K_{Q\C}=Q\K_\C=Q\K$. Let $\hat p:=Q\tp \J p$ and $u\in\Lambda_\K(p)$.
	Since $\lng \hat p-u,\hat p-z\rng=\lng Q\hat p-Qu,Q\hat p-Qz\rng\le 0$,
	for any $z\in\K$, \eqref{eq:L-Proj}
	implies $\Lambda_{\K_{QC}}(p)=\Lambda_{Q\K}(p)=\Lambda_{Q\K}(Q\hat p)=\Lambda_\K(\hat
	p)=\Lambda_K(Q\tp\J p)$.
\end{proof}
Next, we will show how to compute directly the intrinsic $\kappa$-projection for some specific sets using the L-projection. We begin by computing the projection on $\kappa$-hyperbolically convex sets associated with half-spaces.
\begin{example}\label{ex:4}
	Let $r\in\inte(\Lo)$ and $a\in\R^{n+1}$ such that $a\neq 0$ and $a\tp r=0$.
	Consider the hyperplane \[\V:=\big\{x\in\R^{n+1}:a\tp x=0\big\},\] the
	corresponding half-spaces $\V_+:=\lf\{x\in\R^{n+1}:~a\tp x\ge 0\rg\}$, $\V_-:=\lf\{x\in\R^{n+1}:~a\tp x\le 0\rg\}$  and the cone $\K=\V_+\cap\inte(\Lo)$.
	Denote $\bar\K=\V_+\cap\Lo$ the topological closure of $\K$.
	Let $p\in\V_-\cap\inte(\Lo)$.
	By item~\ref{unique} of Corollary \ref{cor:bd},
	$u\in\Lambda_\K(p)\setminus\{0\}$ exists and it  is unique, and  is given by
	\begin{equation}\label{eq:vl}
	u=p-\f{\lng a,\J p\rng}{\lng a,a\rng}\J a.
	\end{equation}
	Let ${{\mathbb H}^n_{\kappa}}$ be the $\kappa$-hyperbolic space form  such that  $\kappa:=-1/{\lng p,p\rng}$. Then \eqref{eq:vl} also provides the $\kappa$-hyperbolic
	projection
	into $\C=\V_+\cap\Hy$ through \eqref{eq:hypp} in  Theorem~\ref{th:hypp}, since
	$\K_\C=\K$, as follows:

		\begin{equation*}
			{P}^{\kappa}_{C}(p)=\frac{\sqrt{\lng a,a\rng}}{
			\sqrt{\kappa\lng a,\J p\rng^2+\lng a,a\rng}}
			\Big(p-\f{\lng a,\J p\rng}{\lng a,a\rng}\J a\Big).
	\end{equation*}
	Indeed, $a\ne 0$ and $a\tp r=0$ with $r\in\inte(\Lo)$ implies $\lng a,a\rng\ne 0$.
	Otherwise, $a\in\Lo\cup-\Lo$. Hence, it follows from Remark \ref{re;PLC}
	that $\lng a,a\rng<0$ or $\lng a,a\rng>0$, which contradicts $\lng
	a,a\rng=0$. On the other hand, $u\in\Lambda(\K)\setminus\{0\}$ if and only if it
	satisfies the
	equation
	\begin{equation}\label{eq:unat}
		u=\Pi_{\bar\K}\lf(\J p+2u_{n+1}e^{n+1}\rg)=\Pi_\V\lf(\J p
		+2u_{n+1}e^{n+1}\rg),
	\end{equation}
	where the first equality of \eqref{eq:unat} comes from item \ref{natural} of
	Corollary \ref{cor:bd}. To justify the second equality of
	\eqref{eq:unat},  note that the first equality of \eqref{eq:unat}
	implies $u_{n+1}\ge 0$. If $u=0$, then $u\in\V\cap\bar\K$ and the second
	equality of \eqref{eq:unat} follows from its first equality because
	$u\in\V\cap\bar\K\subseteq\V$ and $u\in\V\cap\bar\K\subseteq{\bar \K}$ implies
	\[u=\Pi_{\bar\K}\lf(\J p+2u_{n+1}e^{n+1}\rg)=\Pi_{\V\cap\bar\K}\lf(\J p
	+2u_{n+1}e^{n+1}\rg)=\Pi_\V\lf(\J p+2u_{n+1}e^{n+1}\rg).\]
	If $u\ne 0$, then $u_{n+1}>0$ and we obtain  the second equality of \eqref{eq:unat}
	in a similar way as above, because $u\in\V\cap\bar\K\subseteq \bar \K$. Indeed,
	if $u\notin\V$, then $u$ is on $\p\Lo\setminus\bar\K$ which contradicts
	item \ref{bd} of Corollary \ref{cor:bd}. Equation \eqref{eq:unat} together
	with item \ref{natural} of Corollary \ref{cor:bd} implies that
	$u\in\Lambda_\V(p)\setminus\{0\}$. Hence, by using Corollary
	\ref{cor:lph} and item \ref{unique} of Corollary \ref{cor:bd}, we obtain
	formula \eqref{eq:vl}.
\end{example}

In this section, we have conducted a thorough investigation into the computation of the intrinsic $\kappa$-projection into various types of $\kappa$-hyperbolically convex sets. We anticipate that the insights provided herein will be instrumental in advancing and facilitating the practical implementation of numerous methods constrained within this framework, particularly given the historical challenges associated with projection calculations. This exploration naturally leads to considerations of potential applications in the analysis of the gradient projection method, which heavily relies on the intrinsic $\kappa$-projection. Consequently, the subsequent section is dedicated to a comprehensive examination of the gradient projection method, aimed at providing an in-depth analysis of its principles within our established framework.
\section{Gradient projection method on the ${\kappa}$-hyperbolic space forms} \label{sec:socop}
In this section, we focus on the constrained optimization problem on ${\kappa}$-hyperbolic space and present the gradient projection method as a solution approach. To do so, we consider a closed and $\kappa$-hyperbolically convex set ${C} \subset{{\mathbb H}^n_{\kappa}}$ and a differentiable function $f:{{\mathbb H}^n_{\kappa}} \to {\mathbb R}$. The specific form of the constrained optimization problem we are interested in is as follows:
\begin{equation} \label{eq:op}
\quad \displaystyle {\rm Minimize}_{p\in {C}} f(p).
\end{equation}
Let   ${C}^*\neq \varnothing$  be the {\it solution set} of   the problem
\eqref{eq:op} and $-\infty<f^*:= \inf_{p\in {C}}f(p)$   be the  {\it optimum
value}  of $f$.  By using \cite[Proposition 6.1]{FerreiraNemethShu2022} and
equation \eqref{eq:grad}, we
conclude that if  \(\bar p\in {C}\)  is a solution of   \eqref{eq:op}, then
\begin{equation} \label{eq:ncvi}
 \lng \grad f({\bar p}), p \rng = \left\langle {\rm J}f'(p), \Proj_{\bar p}^{\kappa} p \right\rangle \geq 0, \quad \forall ~p\in {C}.
\end{equation}
Any point $\bar p$  satisfying \eqref{eq:ncvi} is called a {\it stationary point} for Problem~\eqref{eq:op}.
In the following corollary   we  present two  important properties of the
projection, which are related to the stationary  points  of the  problem~\eqref{eq:op}.
 \begin{proposition} \label{pr:ProjProperty}
  Let   ${\bar p} \in {C}$ be  such that $\grad f({\bar p})\neq 0$ and  that  $\alpha> 0$.   Then, there hold:
 \begin{enumerate}
	 \item[(i)]\hspace{0mm} The point ${\bar p}$ is stationary for problem~\eqref{eq:op} if and only if ${\bar p}={P}_{C}\left(\exp _{{\bar p}}(-\alpha\grad f({\bar p}))\right)$.
\item[(ii)]\hspace{0mm} If ${p}$ is a nonstationary point for  problem~\eqref{eq:op} then  $\langle \grad f({p}), \log^{\kappa}_{p}{P}_{C}(\exp _{{p}}(-\alpha\grad f({p})))\rangle < 0.$ Equivalently, if there exists    ${\bar \alpha} >0$  satisfying $ \langle \grad f({\bar p}), \log^{\kappa}_{p}{P}_{C}(\exp _{{\bar p}}(-{\bar \alpha}\grad f({\bar p})))\rangle  \geq 0,$ then ${\bar p}$ is stationary for problem~\eqref{eq:op}.
\end{enumerate}
\end{proposition}
\begin{proof}
To prove item~$(i)$,  we first assume that   \(\bar p\in {C}\) is a  stationary  point for Problem~\eqref{eq:op}.  It follows from   \eqref{eq:expinv} and \eqref{eq:ncvi}  that
$
\langle \grad f ({\bar p}), \log^{\kappa}_{\bar p}{P}_{C}(\exp _{{\bar p}}(-\alpha\grad f({\bar p})))\rangle \geq 0.
$
Thus,  by using Proposition~\ref{pr:condw2}, we conclude that  $d({\bar p},{P}_{C}(\exp _{\bar p}(-\alpha\grad f({\bar p}))))\leq 0 $, which implies
that  ${\bar p}={P}_{C}(\exp _{{\bar p}}(-\alpha\grad f({\bar
p})))$.  Conversely,  assume that ${\bar p}={P}_{C}(\exp _{{\bar p}}(-\alpha\grad f({\bar p})))$. Thus, \eqref{d:uexp} implies that  $\langle \log^{\kappa}_{\bar p}\exp _{{\bar p}}\left(-\alpha\grad f({\bar p})\right), \log^{\kappa}_{\bar p}p\rangle \leq 0,$
for all $p\in {C}$, or equivalently, $\langle \alpha\grad f({\bar p}),  \log^{\kappa}_{\bar p}p\rangle \geq 0$, for all $p\in {C}$.  Thus, by using $\alpha
>0$ and \eqref{eq:expinv},  the last inequality implies that  $\langle \grad f (\bar p), \Proj_{\bar p} p \rangle \geq 0,$ for all ~$p\in {C}$. Therefore, the point ${\bar p}$ is stationary for Problem~\eqref{eq:op} and  (i) is proved.  We proceed to prove   item~(ii). Take $p$ a nonstationary point for Problem~\eqref{eq:op}. Thus, item (i) implies that
$
{\bar p}\neq {P}_{C}(\exp _{{\bar p}}(-\alpha\grad f({\bar p}))).
$
Thus, by applying  Proposition~\ref{pr:condw2}, we infer that
 \begin{equation*}
0<\frac{1}{\alpha}d^2\left({\bar p},{P}_{C}\left(\exp _{{\bar p}}\left(-\alpha\grad f({\bar p})\right)\right)\right)\leq - \left\langle \grad f({\bar p}), \log^{\kappa}_{\bar p}{P}_{C}\left(\exp _{{\bar p}}\left(-\alpha\grad f({\bar p})\right)\right)\right\rangle,
\end{equation*}
which implies  the first statement of item~(ii). Finally, note that the second statement of item~(ii) is the contrapositive of the
first statement.
 \end{proof}

\subsection{Gradient  projection method with constant step size}
In this section we analyze the gradient  projection method with constant step size. For that, in this section,  we need to  assume that
\begin{itemize}
\item[{\bf (H1)}] The gradient  vector field $\grad f$  is  Lipschitz continuous  on ${C}$ with   constant $L\geq 0$.
\end{itemize}
 Next  we present the conceptual version of the  {\it gradient projection method} to solve the problem~\eqref{eq:op}.

\begin{algorithm}[H]
\begin{footnotesize}
	\caption{Gradient  projection method on ${{\mathbb H}^n_{\kappa}}$}\label{Alg1s}
	\begin{algorithmic}[1]
		\State Take  $p_0\in {C}$. Set $k=0$;
		\State If $\grad f(p_k)=0$, then {\bf stop} and returns $p_{k}$;  otherwise,   compute a step size  $t_k>0$  such that
		\begin{equation}\label{step.arm}
			p_{k+1}:= {P}^{\kappa}_{C}(\exp^{\kappa}_{p_{k}}\left(-t_{k}\grad f(p_{k})\right)).
		\end{equation}
		\State Update   $k \leftarrow k+1$ and go to   \textbf{2}
	\end{algorithmic}
\end{footnotesize}
\end{algorithm}

It follows from  \eqref{eq:ncvi}  that if  $\grad f(p_k)=0$, then    $p_k$ is a
stationary point for Problem~\eqref{eq:op}.  Therefore, {\it from now on   we
assume that  $\grad f(p_k)\neq0$,   for  all $k=0,1,\ldots$.}  To  proceed, we
assume  that $(p_k)_{k\in {\mathbb N}}$ is  generated by  Algorithm~\ref{Alg1s} with   the following constant stepsize:
\begin{equation} \label{eq:aalpha}
t_k:=\alpha, \qquad \qquad   0<\alpha<\frac{1}{L}, \qquad \quad k=0, 1, \ldots.
\end{equation}
To  simply the  notations of the next results, we define the following positive constant
\begin{equation} \label{eq:CCompP}
\Gamma := (1-\alpha L)/(2\alpha) > 0.
\end{equation}
\begin{lemma} \label{Le:MainConvP}
There hold  $(p_k)_{k\in {\mathbb N}}\subseteq {C}$ and
	\begin{equation} \label{eq:MainIneqP}
		f(p_{k+1}) \leq f(p_{k}) -  \Gamma d_{\kappa}^2(p_{k}, p_{k+1}), \qquad k=0, 1, \ldots.
	\end{equation}
In particular,  the sequence  $(f(p_k))_{k\in {\mathbb N}}$ is non-increasing and converges.
\end{lemma}
\begin{proof}
Since $p_0\in\C$, it follows from  \eqref{step.arm} that  $(p_k)_{k\in {\mathbb
N}}\subseteq\C$. By applying  Lemma~\ref{le:lc} with  $p=p_{k}$ and  $q=p_{k+1}$, we have
$$
f(p_{k+1}) \leq f(p_{k}) + \big\langle\grad f(p_{k}), \log^{\kappa}_{p_k}p_{k+1}\big\rangle +\frac{L}{2}d_{\kappa}^2(p_{k}, p_{k+1}).
$$
Therefore, by using Proposition~\ref{pr:condw2} with $p=p_{k}$ and considering \eqref{step.arm} along with $t_k=\alpha$, we deduce that
\begin{equation*}
		f(p_{k+1}) \leq f(p_{k}) - \frac{1}{2\alpha}d_{\kappa}^2(p_{k}, p_{k+1}) +\frac{L}{2}d_{\kappa}^2(p_{k}, p_{k+1})= f(p_{k}) -  \frac{1-\alpha L}{2\alpha} d_{\kappa}^2(p_{k}, p_{k+1}).
	\end{equation*}
Therefore, we can establish the validity of \eqref{eq:MainIneqP} from
\eqref{eq:CCompP}. Furthermore, the formulas \eqref{eq:aalpha} and
\eqref{eq:MainIneqP} indicate that the sequence $(f(p_k))_{k\in \mathbb{N}}$ is
non-increasing. Additionally, considering that $-\infty < f^*$ and $(f(p_k))_{k\in \mathbb{N}}$ is non-increasing, we can deduce that the sequence converges. Thus, the proof is now complete.
\end{proof}
In the following  we  prove that any  cluster point of $(p_k)_{k\in {\mathbb N}}$  is a solution of the Problem~\eqref{eq:op}.
\begin{theorem}\label{teo.Main}
If ${\bar p}\in {C}$ is a cluster point of the sequence $(p_k)_{k\in {\mathbb N}}$, then  ${\bar {p}}$ is a  stationary point for Problem \eqref{eq:op}.
\end{theorem}
\begin{proof}
If $\grad f(\bar p))=0$, then \eqref{eq:ncvi} implies that  \(\bar p\in {C}\)  is a  stationary point for Problem \eqref{eq:op}.   Now, assume that $\grad
f(\bar p))\neq 0$.   By using  \eqref{eq:MainIneqP}, we obtain that
	\begin{equation} \label{eq:eK1}
		d_{\kappa}^2(p_{k}, p_{k+1}) \leq \frac{1}{\Gamma} \left(f(p_{k}) - f(p_{k+1}) \right), \qquad k = 0,1, \ldots.
	\end{equation}
By applying Lemma~\ref{Le:MainConvP}, we can establish that the sequence
$(f(p_k))_{k\in \mathbb{N}}$ converges. Taking the limit in \eqref{eq:eK1}, we can
deduce that $\lim_{k\to +\infty}d_{\kappa}(p_k, p_{k+1})=0$.
Let ${\bar p}$ be a cluster
point of $(p_k)_{k\in {\mathbb N}}$ and   $(p_{k_j})_{j\in {\mathbb N}}$ a
subsequence of $(p_k)_{k\in {\mathbb N}}$ such that $\lim_{j\to
+\infty}p_{k_j}=~\bar{p}$. As $\lim_{j\to +\infty}d_{\kappa}(p_{k_j+1},
p_{k_j})=0$, we have $\lim_{j\to +\infty}p_{k_j+1}={\bar p}$. On the other hand,
$$
p_{k_j+1} =  {P}_C(\exp^{\kappa}_{p_{k_j}}(-\alpha \grad f(p_{k_j}))), \qquad j=0, 1, \ldots
$$
 and \eqref{d:prjectionef} imply  $\langle  \Proj^{\kappa}_{p_{k_j+1}}\exp^{\kappa}_{p_{k_j}}(-\alpha \grad f(p_{k_j})),   \Proj^{\kappa}_{p_{k_j+1}}q\rangle \leq 0,$ for all  $q\in {C}$. Hence, by setting $v=-\grad f(\bar p)$, using  \eqref{eq:geoexp} and   that ${\bar p}\in {C}\subseteq {{\mathbb H}^n_{\kappa}}$, we have
$
\langle  \Proj^{\kappa}_{\bar p}\exp^{\kappa}_{{\bar p}}(-\alpha \grad f({\bar p})),   \Proj^{\kappa}_{\bar p}q\rangle \leq 0,
$
for all  $q\in {C}$,  which is equivalent  to $\lng \exp^{\kappa}_{\bar p}(-\alpha \grad f(\bar p)),
\Proj^{\kappa}_{\bar p}q \rng \leq 0$, for all $q\in\C$. Therefore, by setting
$v=-\grad f(\bar p)$,  using  \eqref{eq:geoexp} and considering that ${\bar p}\in
{C}\subseteq {{\mathbb H}^n_{\kappa}}$, we can infer that
$$
0\geq \big\langle \cosh(\alpha\sqrt{\kappa}\|v\|){\bar p}+ \sinh(\alpha\sqrt{\kappa}\|v\|) \frac{v}{\sqrt{\kappa}\|v\|},   \Proj^{\kappa}_{\bar p}q \big\rangle= \frac{\sinh(\alpha\sqrt{\kappa}\|v\|)}{\sqrt{\kappa}\|v\|}\lng  v,   \Proj^{\kappa}_{\bar p}q \rng,
$$
for all $q\in\C$. Thus, as  $\sinh(\alpha\sqrt{\kappa}\|v\|)>0$ and $v=-\grad f(\bar p)$,  we have $\langle  \grad f(\bar p),  \Proj^{\kappa}_{\bar p}q \rangle\geq 0,$ for all $q\in\C$,  which,  by using \eqref{eq:ncvi},   implies that  ${\bar p} \in\C$ is a  stationary point for problem~\eqref{eq:op}.
\end{proof}
The item $(i)$ of Proposition~\ref{pr:ProjProperty} implies that if ${p_k}={P}_{C}(\exp _{{p_k}}(-\alpha_k\grad f({p_k})))$, then
$p_k$ is stationary for Problem~\eqref{eq:op}. Since ${p_{k+1}}={P}_{C}(\exp _{{p_k}}(-\alpha_k\grad f({p_k})))$, the distance
between the points $p_k$ and $p_{k+1}$, i.e.,    $d({p_k}, p_{k+1})$, can be seen
as a measure of stationarity of $p_k$.  The next theorem presents an iteration-complexity bound for this measure.
\begin{theorem} 
	For all $N\in \mathbb{N}$ there holds
	$$
	\min\left\{ d_{\kappa}(p_{k}, p_{k+1}):~k=0, 1, \ldots, N \right\}\leq \frac{\sqrt{{f(p_{0})  - f^*}}}{\sqrt{{\Gamma (N+1)}}},
	$$
\end{theorem}
\begin{proof}
	By using \eqref{eq:MainIneqP},   we have $d_{\kappa}^2(p_{k}, p_{k+1})\leq \left(f(p_{k})  -f(p_{k+1}) \right)/\Gamma$,  for all  $k=0, 1, \ldots$.  Since  $f^* \leq f(p_k)$ for all $k$, the last inequality implies
	\begin{equation*}
	\sum_{k=0}^{N} d_{\kappa}^2(p_{k+1}, p_{k})\leq \frac{1}{\Gamma}	\sum_{k=0}^{N} \left(f(p_{k}) - f(p_{k+1})\right)\leq \frac{1}{\Gamma}(f(p_{0})  - f^*) .
	\end{equation*}
	Therefore, $(N+1) \min\left\{d_{\kappa}^2(p_{k}, p_{k+1}):~k=0,1, \ldots, N \right\} \leq (f(p_{0})  - f^*)/\Gamma$, which implies  the desired inequality.
\end{proof}
\subsection{Gradient  projection method with backtracking step size}
In Algorithm~\ref{Alg1s} we need the  Lipschitz constant to compute the stepsize. However,  this  constant  is   not always known or computable.
Therefore, in this section  we  also analyze a version  of Algorithm~\ref{Alg1s}  with a backtracking stepsize rule. The conceptual version of the
gradient projection method on the ${\kappa}$-hyperbolic space   with a backtracking stepsize rule  to solve the problem~\eqref{eq:op} is as follows:

\begin{algorithm}[H]
\begin{footnotesize}
	\caption{ Gradient  projection method on ${{\mathbb H}^n_{\kappa}}$ with backtracking step size rule }\label{AlgsIP}
	\begin{algorithmic}[1]
		\State Take   constants $\rho, \beta, \theta_0\in (0,1)$  and  ${\hat \alpha},{\bar \alpha}>0$. Choose   initial point $p_0\in {C}$.  Set   $k=0$;
		\State If $\grad f(p_k)=0$, then {\bf stop} and returns $p_{k}$; otherwise take
		\begin{equation} \label{eq:coaIP}
			0< {\hat \alpha}\leq {\alpha_{k}}\leq  {\bar \alpha},
		\end{equation}
		\State Compute
		\begin{equation} \label{eq:cpIP}
			y_k:=\exp^{\kappa}_{p_{k}}\left(- {\alpha_{k}}{\grad f(p_{k})}\right), \qquad  \quad z_k:= {P}^{\kappa}_{C}\left(y_k\right),
		\end{equation}
		\begin{equation} \label{eq:jkIP}
			\hspace{-3.5pt}\ell_{k}:=\min \left\{\ell \in {\mathbb N}: ~f\big(q(\beta^{\ell}{\theta_{k}})\big)\leq f(p_{k})+{\rho} \big(\beta^{\ell}{{\theta_{k}} }\big)\big\langle\grad f(p_{k}), \log^{\kappa}_{p_k}z_{k}\big\rangle\right\},
		\end{equation}
		where $q(\tau):= \exp^{\kappa}_{p_{k}}\big(\tau \log^{\kappa}_{p_k}z_k\big)$ denotes  the geodesic segment joining $p_k$ to $z_k$;
		\State Set   $p_{k+1}:= q(\beta^{\ell_k}{\theta_{k}})$ and  $\theta_{k+1}:= \beta^{\ell_{k}-1}{\theta_{k}}$;
		\State Update   $k \leftarrow k+1$ and go to {\bf 2}.
	\end{algorithmic}
	\end{footnotesize}
\end{algorithm}

In the next proposition we prove that  Algorithm~\ref{AlgsIP}  is well defined
and two useful inequalities.
\begin{proposition}  \label{pr:wdAstsIP}
The Algorithm~\ref{AlgsIP} is well defined  and  generates a sequence
$(p_k)_{k\in \mathbb N}\subseteq C$. Moreover, the following inequality holds
\begin{equation} \label{eq:ddIP}
\big\langle \grad f(p_k),  \log^{\kappa}_{p_k}z_k\big\rangle\leq -\frac{1}{\alpha_k}d^2(p_k, z_k)<0,  \qquad \forall k\in {\mathbb N}.
\end{equation}
Furthermore,
\begin{equation} \label{eq:inqmIP}
f(p_{k+1})\leq f(p_{k})+{\rho} \beta {\theta_{k+1}} \big\langle\grad f(p_{k}),\log^{\kappa}_{p_k}z_{k}\big\rangle, \qquad \forall k\in {\mathbb N}.
\end{equation}
\end{proposition}
\begin{proof}
Bearing in mind that  $p_0\in {C}$, without loss of generality  we can
assume that $p_k\in {C}$.  Since  ${C} \subset{{\mathbb H}^n_{\kappa}}$
is   closed, ${P}_{C}\left(y_k\right)$ is a singleton.  Hence,   the
point $z_k$ in  \eqref{eq:cpIP}  is well defined and  belongs to the set  ${C}$. Since $y_k$ are  given by \eqref{eq:cpIP}, Proposition~\ref{pr:condw2}  with
$v=\grad f(p_k)$,  $p=p_k$ and $\alpha= {\alpha_{k}}$ implies that
$$
\big\langle \grad f(p_k),  \log^{\kappa}_{p_k} {P}_{C}\left(y_k\right)\big\rangle\leq -\frac{1}{\alpha_k}d_{\kappa}^2(p_k, {P}_{C}\left(y_k\right)).
$$
Thus, due to  $z_k= {P}_{C}\left(y_k\right)$, the inequality
\eqref{eq:ddIP} follows from the previous one.  From the differentiability of $f$ is and
$0<\beta<1$, we conclude that
\begin{equation} \label{eq:da3fc2IP}
\lim_{\ell \to +\infty}\frac{f(\exp^{\kappa}_{p_{k}}\big(\beta^{\ell}{\theta_{k}} \log^{\kappa}_{p_k}z_k\big))- f(p_{k})}{\beta^{\ell}{{\theta_{k}} }}=\big\langle \grad f(p_k),  \log^{\kappa}_{p_k}z_k\big\rangle.
\end{equation}
Thus, taking into account that $0<\rho<1$ and using item (ii) of Proposition~\ref{pr:ProjProperty}, the combination of \eqref{eq:ddIP} with \eqref{eq:da3fc2IP} implies that
\begin{equation*}
\lim_{\ell \to +\infty}\frac{f(\exp^{\kappa}_{p_{k}}\big(\beta^{\ell}{\theta_{k}} \log^{\kappa}_{p_k}z_k\big))- f(p_{k})}{\beta^{\ell}{{\theta_{k}} }}< {\rho}\big\langle \grad f(p_k),  \log^{\kappa}_{p_k}z_k\big\rangle.
\end{equation*}
Therefore,  there exists ${\hat \ell}\in  {\mathbb N}$  such that
$$
f(\exp^{\kappa}_{p_{k}}\big(\beta^{\ell}{\theta_{k}} \log^{\kappa}_{p_k}z_k\big))< f(p_{k})+{\rho} \big(\beta^{\ell}{{\theta_{k}} }\big)\big\langle \grad f(p_k),  \log^{\kappa}_{p_k}z_k\big\rangle,  \qquad \forall \ell\in   {\mathbb N};~   {\hat \ell}\leq \ell,
$$
which implies $\ell_k$ at {\bf Step 3} is well defined,  and then
Algorithm~\ref{AlgsIP}  is also well defined.  The definition of the next iterate
$p_{k+1}$ implies that it belongs to the geodesic through $p_k$ and $z_k$. Since
$p_k$  and  $z_k$  belong to ${C}$  and ${C}$ is $\kappa$-hyperbolically convex, we infer that $p_{k+1}$ also belongs to ${C}$. Therefore, the  generated  sequence $(p_k)_{k\in \mathbb N}$ belongs to ${C}$. The inequality \eqref{eq:inqmIP} follows from  the definition of $\ell_k$ at {\bf Step 3} and definitions of $p_{k+1}$ and $\theta_{k+1}$ at {\bf Step 4}, which concludes  the proof.
\end{proof}
\begin{theorem}\label{teo.MainIP}
If ${\bar p}\in {C}$ is a cluster point of the sequence $(p_k)_{k\in {\mathbb N}}$, then  ${\bar {p}}$ is a  stationary point for Problem \eqref{eq:op}.
\end{theorem}
\begin{proof}
 Let ${\bar p}$ be a cluster point of the sequence  $(p_k)_{k\in {\mathbb N}}$. If $\grad f(\bar p)=0$, then \eqref{eq:ncvi} implies that  \(\bar p\in {C}\)  is a  stationary point for Problem \eqref{eq:op}.   Now, assume that $\grad f(\bar p)\neq 0$. Combining     \eqref{eq:ddIP} with \eqref{eq:inqmIP}   we have
 \begin{equation} \label{eq:inqmFIPnl}
0<-{\rho} \beta {\theta_{k+1}} \big\langle\grad f(p_{k}),\log^{\kappa}_{p_k}z_{k}\big\rangle \leq f(p_{k})-f(p_{k+1}), \qquad \forall k\in {\mathbb N}.
\end{equation}
It follows from  \eqref{eq:inqmFIPnl}  that  $(f(p_k))_{k\in {\mathbb N}}$ is non-increasing and  considering that $-\infty <f^*$, we conclude that it converges.
Since $(f(p_k))_{k\in {\mathbb N}}$ converges,  we obtain   from \eqref{eq:inqmFIPnl} that
\begin{equation} \label{eq:fcddIP}
\lim_{k\to +\infty}  {\theta_{k+1}}\big\langle\grad f(p_{k}), \log^{\kappa}_{p_k}z_{k}\big\rangle=0.
\end{equation}
Since ${\bar p}$ is a cluster point of  sequence   $(p_k)_{k\in {\mathbb N}}$, we take $(p_{k_j})_{j\in {\mathbb N}}$   a subsequence of
$(p_k)_{k\in {\mathbb N}}$ such that $\lim_{j\to +\infty}p_{k_j}=~\bar{p}$. Considering that   $(\alpha_k)_{k\in {\mathbb N}}\subseteq [{\hat
\alpha}, {\bar \alpha}]$ and $(\theta_k)_{k\in {\mathbb N}}\subseteq (0, 1)$,  we can assume without generality that
\begin{equation} \label{eq:alphthet}
\lim_{j\to +\infty}\alpha_{k_j}=:{\tilde \alpha}\geq {\hat \alpha}>0, \qquad \lim_{j\to +\infty}\theta_{k_j+1}=:{\tilde \theta}\in [0, 1].
\end{equation}
 It follows from \eqref{eq:cpIP} that  $z_{k_j}:= {P}^{\kappa}_{C}\big(\exp^{\kappa}_{p_{k_j}}(- {\alpha_{k_j}}{\grad f(p_{k_j})})\big)$. Thus,
 by using Proposition~\ref{pr:PrjCont}, taking into account that $\lim_{j\to
 +\infty}p_{k_j}=~\bar{p}$,   and \eqref{eq:alphthet}  we infer that
\begin{equation} \label{eq:cpIPnl}
\lim_{j\to +\infty} z_{k_j}:= {P}^{\kappa}_{C}\big(\exp^{\kappa}_{\bar{p}}(- {\tilde \alpha}{\grad f(\bar{p})})\big)=:{\bar z}.
\end{equation}
From \eqref{eq:fcddIP} we obtain  that $\lim_{k\to +\infty}
{\theta_{k_j+1}}\big\langle\grad f(p_{k_j}),
\log^{\kappa}_{p_{k_j}}z_{k_j}\big\rangle=0$. Hence, the combination of \eqref{eq:alphthet}  and  \eqref{eq:cpIPnl}  with  $\lim_{j\to +\infty}p_{k_j}=~\bar{p}$  yields
\begin{equation} \label{eq:Mainnl}
{\tilde \theta}\big\langle\grad f({\bar p}), \log^{\kappa}_{\bar p}{\bar z}\big\rangle=0.
\end{equation}
Since ${\tilde \alpha}\in [0, 1]$,  we have two possibilities that we must consider: ${\tilde \theta}>0$ or ${\tilde \theta}=0$. First we assume that ${\tilde \theta}>0$. In this case, by using \eqref{eq:ddIP} we have
\begin{equation} \label{eq:pddIPnl}
d^2(p_{k_j}, z_{k_j})\leq -{\alpha_{k_j}}\big\langle \grad f(p_{k_j}),  \log^{\kappa}_{p_{k_j}}{z_j}\big\rangle,  \qquad \forall j\in {\mathbb N}.
\end{equation}
By taking limit in   \eqref{eq:pddIPnl}, and by using \eqref{eq:alphthet}
and  \eqref{eq:cpIPnl},   together with 
\eqref{eq:Mainnl}, we get  $\lim_{k\to +\infty} d(p_{k_j}, z_{k_j})=0$. Thus,  due to  $\lim_{j\to +\infty}p_{k_j}=~\bar{p}$ and  $\lim_{j\to +\infty}z_{k_j}={\bar z}$, we have ${\bar z}={\bar p}$.  Hence, using \eqref{eq:cpIPnl} we obtain that
$$
{\bar p} = {P}^{\kappa}_{C}\big(\exp^{\kappa}_{\bar{p}}(- {\tilde \alpha}{\grad f(\bar{p})})\big).
$$
Therefore, item $(i)$ of Proposition~\ref{pr:ProjProperty} implies that  ${\bar {p}}$ is a  stationary point for Problem \eqref{eq:op}.

Now,  we assume that ${\tilde \theta}=0$. In this case, we consider the following auxiliary sequence $(q_{k_j})_{j\in {\mathbb N}}$ defined by
\begin{equation} \label{eq:auxnl}
q_{k_j}:= \exp^{\kappa}_{p_{k_j}}\big({\theta_{k_j+1}} \log^{\kappa}_{p_{k_j}}z_{k_j}\big),  \qquad \forall k\in {\mathbb N},
\end{equation}
Since  $\lim_{j\to +\infty}\theta_{k_j}=0$, we have $\lim_{j\to +\infty}q_{k_j}={\bar p}$. Moreover,  the definition of $\ell_{k}$ in \eqref{eq:jkIP} implies that
$
f(q_{k_j})> f(p_{k_j})+{\rho} {\theta_{k_j+1}} \big\langle\grad f(p_{k_j}), \log^{\kappa}_{p_{k_j}}z_{k_j},\big\rangle.
$
Via combining \eqref{eq:auxnl}  with the last inequality, we obtain that
\begin{equation*}
\frac{f\big( \exp^{\kappa}_{p_{k_j}}({\theta_{k_j+1}} \log^{\kappa}_{p_{k_j}}z_{k_j})\big)-f(p_{k_j})}{\theta_{k_j+1}}> {\rho} \big\langle\grad
f(p_{k_j}), \log^{\kappa}_{p_{k_j}}z_{k_j},\big\rangle.
\end{equation*}
Since $\lim_{j\to +\infty}\theta_{k_j+1}=0$, $\lim_{j\to
+\infty}p_{k_j}=~\bar{p}$, $\lim_{j\to +\infty}z_{k_j}={\bar z}$ and $f$ is
differentiable, By taking limit in the last inequality we arrive to
$
\big\langle\grad f(\bar{p}), \log^{\kappa}_{\bar{p}}\bar{z},\big\rangle\geq {\rho} \big\langle\grad f(\bar{p}), \log^{\kappa}_{\bar{p}}\bar{z},\big\rangle.
$
Due to $0<\rho<1$, it follows from the last inequality that $\big\langle\grad f(\bar{p}), \log^{\kappa}_{\bar{p}}\bar{z},\big\rangle\geq 0$. Thus, using \eqref{eq:cpIPnl}  we obtain  that
$$
\big\langle\grad f(\bar{p}), \log^{\kappa}_{\bar{p}}{P}^{\kappa}_{C}\big(\exp^{\kappa}_{\bar{p}}(- {\tilde \alpha}{\grad f(\bar{p})})\big)\big\rangle\geq 0.
$$
Therefore, using the  inequality in item (ii) of Proposition~\ref{pr:ProjProperty} we deduce that  ${\bar {p}}$ is a  stationary
point for Problem \eqref{eq:op}. Thus, the proof is now complete.
\end{proof}
To proceed, we need  the assumption  (H1). For simplifying  the notation, it is convenient to define the following positive constant
\begin{equation} \label{eq:baraphaIP}
 {\bar \theta}:=\min \Big\{ \theta_{0},\frac{2}{{\bar \alpha}L}(1-\rho)\Big\}.
\end{equation}
\begin{lemma} \label{Le:MainConvPAstsIP}
 The following  inequality holds:
\begin{equation} \label{eq:baphaIP}
\theta _{k}\geq {\bar \theta}, \quad \forall k \in \mathbb{N}.
\end{equation}
\end{lemma}
\begin{proof}
Since $p_0\in {C}$, it follows from Proposition~\ref{pr:wdAstsIP} that
$(p_k)_{k\in {\mathbb N}}\subseteq {C}$.  The  inequality \eqref{eq:baphaIP}
immediately holds for $k=0$.   Lemma~\ref{le:lc} with ${D}={{\mathbb
H}^n_{\kappa}}$, $p=p_{k}$ and  $q=\exp _{p_{k}}(\beta^{\ell_k-1}{\theta_{k}}
\log^{\kappa}_{p_k}z_k)$ yield
$$
f(\exp^{\kappa} _{p_{k}}\big(\beta^{\ell_k-1}{\theta_{k}} \log^{\kappa}_{p_k}z_k\big)) \leq f(p_{k}) + \beta^{\ell_k-1}{\theta_{k}}\langle\grad f(p_{k}), \log^{\kappa}_{p_k}z_{k}\rangle +\frac{L}{2}(\beta^{\ell_k-1}{\theta_{k}})^2d^2(p_{k}, z_{k}).
$$
On the other hand, it follows from the definition of   $\ell_k$ in \eqref{eq:jkIP} that
$$
f(\exp^{\kappa} _{p_{k}}\big(\beta^{\ell_k-1}{\theta_{k}} \log^{\kappa}_{p_k}z_k\big)) > f(p_{k}) + \rho\beta^{\ell_k-1}{\theta_{k}}\langle\grad f(p_{k}), \log^{\kappa}_{p_k}z_{k}\rangle.
$$
Hence, by combining the two previous inequalities and  keeping in mind that
$\theta_{k+1}= \beta^{\ell_{k}-1}{\theta_{k}}$,  we get
$$
 -(1-\rho)\theta_{k+1}\langle\grad f(p_{k}), \log^{\kappa}_{p_k}z_{k}\rangle <\frac{L}{2}\theta_{k+1}^2d^2(p_{k}, z_{k}).
$$
Having in mind that $0<\rho<1$, by  applying Proposition~\ref{pr:wdAstsIP},  we
infer from the last inequality  that
$$
(1-\rho)\theta_{k+1}\frac{1}{\alpha_k}d^2(p_{k}, z_{k})<\frac{L}{2}\theta_{k+1}^2d^2(p_{k}, z_{k}).
$$
Thus, due to  $0<{\alpha_{k}}< {\bar \alpha}$, $\theta_{k+1}\neq 0$ and $d(p_{k}, z_{k})\neq 0$,   it follows from the previos inequality that
$
({2}/{{\bar \alpha}L})(1-\rho)<\theta_{k+1}.
$
Therefore, by using   \eqref{eq:baraphaIP} we have  $ {\bar \theta}<\theta_{k+1}$. Because   \eqref{eq:baphaIP}  holds for $k=0$, the proof is complete.
\end{proof}
As the quantity $d(p_k, z_k)$ can be interpreted as a measure of the stationarity of $p_k$, the following theorem provides a bound on the iteration complexity for this measure.
\begin{theorem} 
	For all $N\in \mathbb{N}$ there holds
	$$
	\min\left\{ d({p_k}, {z_k}):~k=0, 1, \ldots, N \right\}\leq \frac{\sqrt{{\bar \alpha}(f(p_{0})  - f^*)}}{\sqrt{{\rho} \beta {\bar \theta}}}   \frac{1}{ \sqrt{(N+1)}}.
	$$
\end{theorem}
\begin{proof}
	Inequalities  \eqref{eq:inqmIP},  \eqref{eq:coaIP}, \eqref{eq:ddIP}  and
	\eqref{eq:baphaIP} yield

\begin{equation*}
{\rho} \beta {\bar \theta} \frac{1}{\bar \alpha}d^2(p_k, z_k)\leq -{\rho} \beta {\theta_{k+1}} \big\langle \grad f(p_k),  \log^{\kappa}_{p_k}z_k\big\rangle \leq f(p_{k})-f(p_{k+1}), \qquad \forall k\in {\mathbb N}.
\end{equation*}
Since  $f^* \leq f(p_k)$ for all $k$, the last inequality implies
	\begin{equation*}
	\sum_{k=0}^{N} d^2({p_k}, {z_k})\leq 	 \frac{\bar \alpha}{{\rho} \beta {\bar \theta}}\sum_{k=0}^{N} \left(f(p_{k}) - f(p_{k+1})\right)\leq \frac{\bar \alpha}{{\rho} \beta {\bar \theta}}(f(p_{0})  - f^*) .
	\end{equation*}
	Therefore, $(N+1) \min\left\{d^2({p_k}, {z_k}):~k=0,1, \ldots, N \right\} \leq {\bar \alpha}(f(p_{0})  - f^*)/({{\rho} \beta {\bar \theta}})$, which implies  the desired inequality.
\end{proof}
\section{The constrained center of mass in $\kappa$-hyperbolic space forms} \label{sec:App}
Let us consider the problem of computing the center of mass, which is also known
as the Fermat-Weber problem, the Fr\'echet mean, or the variance problem (see \cite{AfsariVidal2013, Bacak2014}). This problem arises frequently in various practical applications without constraints, as demonstrated in \cite{Afsari2011, fletcher2004principal, pmlr-v119-lou20a}. However, in this section, we will also consider the constrained version of this problem, which may become important in the future.

Let $\{q_1, \dots, q_m\}\subseteq {{\mathbb H}^n_{\kappa}}$ be a data set of distinct points,  $q_i\neq q_j$ for $i\neq j$,
and  $\mu_1, \ldots \mu_m$ be positive  weights with $\sum_{j=1}^m\mu_j=1$ and $\sigma>1$. Let ${C}\subseteq  {{\mathbb H}^n_{\kappa}}$ be a closed and $\kappa$-hyperbolically  convex set.  The constrained \emph{ center of mass}  minimization problem is stated as follows:
\begin{equation}\label{eq:probminfv}
 \arg\min_{p\in{ C}} \sum_{i=1}^m \mu_i d_{\kappa}^{\sigma}(p,q_i),
\end{equation}
where $d_{\kappa}$ is the  intrinsic distance on the  $\kappa$-hyperbolic space
form defined in \eqref{eq:Intdist}.  Assume that ${\sigma}\geq 2$.  Then, the
objective function of problem~\eqref{eq:probminfv} is continuously
differentiable.  For simplifying the notations, let $ \zeta: {{\mathbb H}^n_{\kappa}} \to {\mathbb R}$  and $ \zeta_i:{{\mathbb H}^n_{\kappa}} \to {\mathbb R}$ be defined, respectively,  by
\begin{equation} \label{eq:fcm}
 \zeta(p):=\sum_{i=1}^m \mu_i \zeta_i (p), \qquad \quad  \zeta_i (p):=d_{\kappa}^{\sigma}(p,q_i),
\end{equation}
for all $p\in  {{\mathbb H}^n_{\kappa}}$ and $i=1, \ldots, m$. The  sub-level set
of $\zeta$ associated to $c>0$ is defined by
\begin{equation*}
L_{\zeta,{ C}}(c):=\{p\in { C}:~  \zeta(p)\leq c\}.
\end{equation*}
\begin{proposition}  \label{eq:ExistSol}
For any $c>0$, the sub-level $L_{\zeta,{ C}}(c)$ is a compact set. As a consequence, the problem \eqref{eq:probminfv} has a global  solution.
\end{proposition}
\begin{proof}
We can assume without loss of generality that  the sub-level set
$L_{\zeta,{ C}}(c)$ is a nonempty set.  Let $0<\mu_0:=\min\{\mu_1, \ldots,
\mu_m\}$, $r:=(c/\mu_0)^{1/\sigma}$ and  $B(q_i, r):=\{p\in {{\mathbb
H}^n_{\kappa}}:~ d_{\kappa}(p,q_i)\leq r\}$, for $i=1, \ldots, m$.  Hence, we
have $L_{\zeta,{ C}}(c) \subseteq \cup_{i=1}^{m}B(q_i, r)$. Consequently,
$L_{\zeta,{ C}}(c)$ is bounded. Since $\zeta$ is continuous, we conclude
that $L_{\zeta,{ C}}(c)$ is also closed. Therefore, $L_{\zeta,{
C}}(c)$ is compact, which proves the first part.  For the second part,
note that since the sub-level $L_{\zeta,{ C}}(c)$ set is compact and nonempty,
the constrained optimization problem $\arg\min_{p\in{ L_{\zeta,{ C}}(c)}}\zeta(p)$
has a global solution. Let us denote this solution by ${\bar p}\in L_{\zeta,{
C}}(c)$. Thus, we have  $\zeta({\bar p})\leq  \zeta(p)\le c$ if       $p\in L_{\zeta,{
C}}(c)$ and,   for all  $p\in { C}\backslash  L_{\zeta,{ C}}(c)$ we have
$\zeta(p)>c\geq \zeta({\bar p})$. Therefore,   $\zeta({\bar p})\leq  \zeta(p)$  for all $p\in { C}$, and the problem \eqref{eq:probminfv} has a global  solution.
\end{proof}
For proving  the next proposition we will need the following definition.
\begin{definition}\label{def:cf-b}
Let ${C}\subseteq {{\mathbb H}^n}$ be a $\kappa$-hyperbolically convex set and $I\subseteq {\mathbb R}$ an interval.
A function $f:{C}\to {\mathbb R}$  is called $\kappa$-hyperbolically convex
(respectively, strictly $\kappa$-hyperbolically convex) if for any   geodesic segment $\gamma:I\to {C}$, the composition
$ f\circ \gamma :I\to {\mathbb R}$ is convex (respectively, strictly convex) in the usual sense.
\end{definition}
\begin{proposition}  \label{eq:Uniq}
The problem \eqref{eq:probminfv} has  global solution and the solution is unique.
\end{proposition}
\begin{proof}
Let  $p_0\in{C}$. Thus,  the sublevel set $\{p\in  {{\mathbb H}^n_{\kappa}}:~  \zeta(p)\leq  \zeta(p_0)\}\neq \varnothing$. Thus, Proposition~\ref{eq:ExistSol} implies that the problem \eqref{eq:probminfv} has global solution. We proceed to prove the uniqueness.
First note that, due to ${\sigma}\geq 2$, \eqref{eq:Hess} and similar
arguments  to the ones used to prove   \cite[Lemma 2.4]{FerreiraNemethShu2022}, we
can prove that the Hessian $\Hess \zeta_i(p)$ is positive deﬁnite. Thus
\cite[Proposition~5.4]{FerreiraNemethShu2022} implies that $\zeta_i$  is strictly
$\kappa$-hyperbolically convex.
Therefore, as $\mu_1, \ldots, \mu_m$ are positive, it can be inferred from \eqref{eq:fcm} that $\zeta$ is strictly $\kappa$-hyperbolically convex.
As a consequence, the global solution of  problem~\ref{eq:probminfv} is unique.
\end{proof}
Now, let us analyze the convergence of Algorithm~\ref{AlgsIP} when it is applied to the problem \eqref{eq:probminfv}.
\begin{theorem}\label{teo.MainCM}
Consider the sequence $(p_k)_{k\in \mathbb{N}}$ generated by Algorithm~\ref{AlgsIP} to solve problem \eqref{eq:probminfv}.
Then,  $(p_k)_{k\in {\mathbb N}}$ converges to the solution of problem \eqref{eq:probminfv}.
\end{theorem}
\begin{proof}
Consider $p_0$ as the starting point of the sequence $(p_k)_{k\in \mathbb{N}}$.
From Proposition~\ref{pr:wdAstsIP}, we can deduce that the sequence $(\zeta(p_k))_{k\in \mathbb{N}}$ is non-increasing. Hence, we have
$
	(p_k)_{k\in \mathbb{N}} \subseteq \{p\in \mathbb{H}^n_{\kappa}:~ \zeta(p) \leq \zeta(p_0)\}.
$
By applying the first part of Proposition~\ref{eq:ExistSol}, we can establish
that $(p_k)_{k\in \mathbb{N}}$ is bounded. Thus, $(p_k)_{k\in \mathbb{N}}$ has a cluster
point $\bar{p}\in {C}$. Utilizing Theorem~\ref{teo.MainIP}, we conclude that
$\bar{p}$ is a stationary point for the problem~\eqref{eq:probminfv}. Since $\zeta$ is $\kappa$-hyperbolically convex
and stationary points are global solutions, $\bar{p}$ is also a solution for the
problem~\eqref{eq:probminfv}. Therefore, since Proposition~\ref{eq:Uniq} implies that problem
\eqref{eq:probminfv} has only one solution, the sequence $(p_k)_{k\in \mathbb{N}}$ has only
one cluster point and must converge.
\end{proof}

Let us proceed to provide a formula for the gradient of the objective function of problem
\eqref{eq:probminfv}. For $p\neq q_i$, \eqref{eq:Intdist} and \eqref{eq:grad}
imply
\begin{equation*}
\grad  \zeta_i(p)=-\sigma d^{\sigma -2}_{\kappa}(p,q_i)\log^{\kappa}_{p}{q_i}, \qquad \quad  i=1, \ldots, m.
\end{equation*}
Thus,  due to ${\sigma}\geq 2$,  the objective function of
problem~\eqref{eq:probminfv} is continuously differentiable.  Hence, it follows,
from the last equality, \eqref{eq: proj}, \eqref{eq:Intdist}  and
\eqref{eq:expinv}, that   the gradient of the objective function in \eqref{eq:probminfv} is given by
 \begin{equation} \label{eq:gradfcm}
  \grad   \zeta(p)=  -\sum_{i=1}^m \mu_i  \frac{\sigma{\kappa}^{1-\frac{\sigma}{2}}({\rm arcosh} (-\kappa \lng p , q_i\rng))^{\sigma-1}}{\sqrt{{\kappa}^2\lng {p}, {q_i}\rng^2-1}}  \left( {\rm I}+{\kappa}_{p}{p}^T{\rm J}\right) {q_i}.
\end{equation}
To implement Algorithm~\ref{AlgsIP}, it is necessary to utilize the formula to project  intrinsically  into the $\kappa$-hyperbolically convex set  ${C}$ in consideration, some of which were provided  in Section~\ref{sec:pj},   and \eqref{eq:gradfcm}.

\section{Numerical Experiments} \label{sec:Nuerics}

As a first example we consider the Riemannian Center of mass~\cite{Karcher:1977}
constrained to a subset:
For some $p \in \mathbb H_{\kappa}^n$ and a radius $r  > 0$ we denote by
\begin{equation*}
	C_{p,r} := \{ q \in \mathbb H_{\kappa}^n: d_{\kappa}(p, q) \leq r \}
\end{equation*}
the ball around $p$ of radius $r$.

Let $q_i \in \mathbb H^2$, $i=1,\ldots,N$, $N \in \mathbb N$,
be a given set of data points. Then the constrained (weighted) Riemannian center of
mass is given by
\begin{equation}
	\arg\min_{p\in{ C}} f(p) \qquad\text{ with }\qquad f(p) = \sum_{i=1}^m \mu_i d_{\kappa}^2(p,q_i).
\end{equation}

We compare the Augmented Lagrangian Method (ALM) and the Exact Penalty Method (EPM)~\cite{LiuBoumal2020}
to our projected gradient (PGA) algorithm.
All three algorithms are implemented using the programming language Julia in
\texttt{Manopt.jl}~\cite{Bergmann:2022}\footnote{Version 0.5.12, the documentation of the projected gradient method is available at \href{https://manoptjl.org/stable/solvers/projected_gradient_method/}{manoptjl.org/stable/solvers/projected\_gradient\_method/}}
and the hyperbolic manifold in \texttt{Manifolds.jl}\footnote{Version 0.10.14, see \url{https://juliamanifolds.github.io/Manifolds.jl/stable/}.}~\cite{AxenBaranBergmannRzecki:2023}.
All time measurements were done using \texttt{Julia} version 1.11.4, on a MacBook Pro with an Apple M1 chip (2021), 16 GB Ram,
running Mac OS Sequoia 15.3.2.

For the first experiment\footnote{
The complete code is available at \href{https://juliamanifolds.github.io/ManoptExamples.jl/stable/examples/Constrained-Mean-H2/}{juliamanifolds.github.io/ManoptExamples.jl/stable/examples/Constrained-Mean-H2/}
} we consider $\kappa = 1$, $d=2$, $r=1$, and $N=200$ data points
and a constrained set of $C = C_{p_0,r}$ woth $p_0 = (0,0,1)^{\mathrm{T}}$.
These are created as $X_i = (\frac{3}{2} + \sigma x_1, \frac{3}{2} + \sigma x_2, 0)^{\mathrm{T}} \in T_{p_0}\mathbb H_1^2$,
$i=1,\ldots,N$,
where $\sigma = \frac{3}{2}$ and $x_1, x_2$ are random normally distributed independent values with mean $0$ and variance $1$.
We then set $q_i = \exp_{p_0} X_i$.

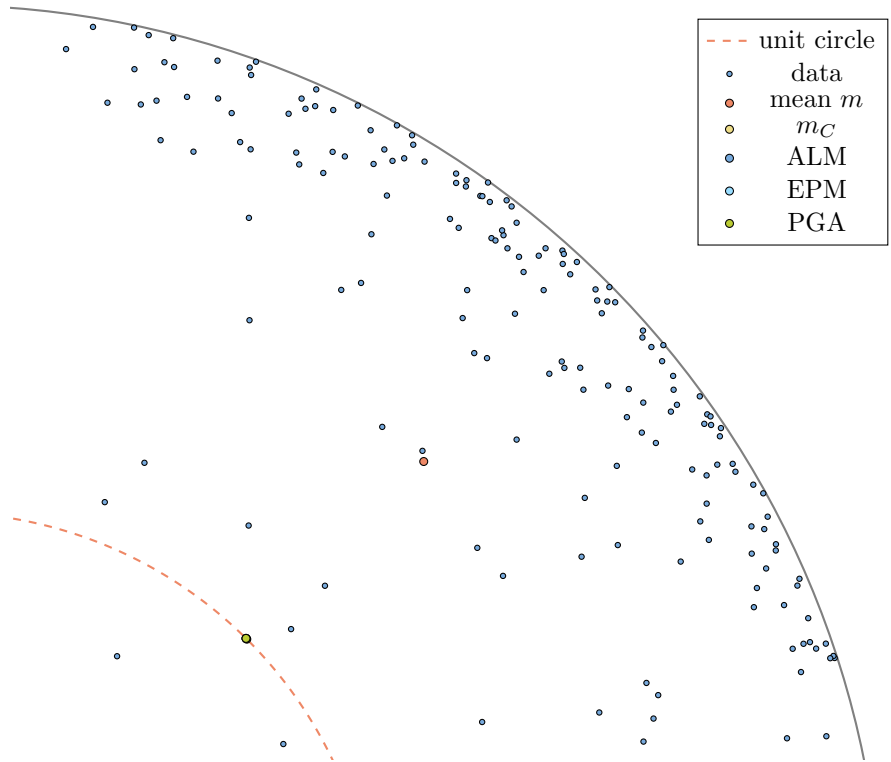
\begin{figure}
    \pgfplotstableread[col sep = comma]{C2Dball-mean-200-pts.csv}\CTwoDBalldata
    \pgfplotstableread[col sep = comma]{C2Dball-mean-200-constraint.csv}\CTwoDBallConstr
    \pgfplotstableread[col sep = comma]{C2Dball-mean-200-epm.csv}\CTwoDBallEPM
    \pgfplotstableread[col sep = comma]{C2Dball-mean-200-mean.csv}\CTwoDBallM
    \pgfplotstableread[col sep = comma]{C2Dball-mean-200-pg.csv}\CTwoDBallPG
    \pgfplotstableread[col sep = comma]{C2Dball-mean-200-proj.csv}\CTwoDBallPM
    \pgfplotstableread[col sep = comma]{C2Dball-mean-200-alm.csv}\CTwoDBallALM
    \centering
    \begin{tikzpicture}
		\begin{axis}[
		axis equal,  
		width = .8\textwidth,  
		xmin = 0.1, xmax = 1.01,
		ymin = 0.2, ymax = 1.01,
		axis lines = none, 
		legend pos=north east
		]
		\draw[gray, thick] (axis cs: 0, 0) circle [radius = 100];
	    \addplot [thick,TolLightOrange, dashed] table [x=x, y=y, col sep=comma] {\CTwoDBallConstr};
		\addlegendentry{unit circle}
	    \addplot [only marks,mark=*,mark options ={fill=TolLightBlue, draw=none, scale=0.5}] table [x=x, y=y, col sep=comma] {\CTwoDBalldata};
           \addlegendentry{data}
		\addplot [only marks,mark=*,mark options ={fill=TolLightOrange, draw=none, scale=0.75}] table [x=x, y=y, col sep=comma] {\CTwoDBallM};
        \addlegendentry{mean $m$}
		\addplot [only marks,mark=*,mark options ={fill=TolLightYellow, draw=none, scale=0.75}] table [x=x, y=y, col sep=comma] {\CTwoDBallPM};
        \addlegendentry{$m_C$}
		\addplot [only marks,mark=*,mark options ={fill=TolLightBlue, draw=none, scale=0.75}] table [x=x, y=y, col sep=comma] {\CTwoDBallALM};
        \addlegendentry{ALM}
		\addplot [only marks,mark=*,mark options ={fill=TolLightCyan, draw=none, scale=0.75}] table [x=x, y=y, col sep=comma] {\CTwoDBallEPM};
        \addlegendentry{EPM}
		\addplot [only marks,mark=*,mark options ={fill=TolLightPear, draw=none, scale=0.75}] table [x=x, y=y, col sep=comma] {\CTwoDBallPG};
        \addlegendentry{PGA}
        \end{axis}
    \end{tikzpicture}
    \caption{Illustration a section of the Poincaré ball, the majority of the sampled points outside the constrained set illustrated as an orange dashed circle.
	The Riemannian center of mass lies outside the constrained set. The projected mean and all minimisers lie very close together, while the projected mean
	is still higher in cost and a reasonable distance away.}
	\label{fig:H2Cmean:illustration}
\end{figure}

For the $\mathbb H_1^2$ we can use the representation in the Poincaré ball to visualise the data,
which are shown in~Figure\ref{fig:H2Cmean:illustration}. The mean $m \not\in C$ is projected onto the set
$C$ to obtain $m_C$.

For the two comparison algorithms running constrained optimisation we formulate the constrained in a smooth way as

\begin{eqnarray}
	g(p) = d_1^2(p_0,p) - r^2 \quad\text{ and hence }\quad \operatorname{grad} g(p) = -2\log_{p}^1p_0.
\end{eqnarray}

and we use the default stopping criteria for both ALM and EPM. $p_0$ also serves as a start point for all three algorithms

For PGA we use
a constant stepsize $\alpha_k = 1$, an Armijo-type backtracking for the step size $\ell_k$ in \eqref{eq:jkIP},
and a stopping criterion of having reached a stationary point, that is when $d_\kappa(p_{k-1}, z_k) \leq \epsilon$ with $\varepsilon = 10^{-7}$.

The resulting points of all three algorithm runs are also displayed in Figure\ref{fig:H2Cmean:illustration}
as ALM, EPM, and PGA, respectively.
These are visually indistinguishable also from the projected mean.
Between the 3 results the maximal distance is $1.5\cdot10^{-5}$,
ALM and PGA are even as close as $9.36\cdot10^{-8}$.
Their distance to the projected mean is at least $2.613\cdot10^{-3}$, so the result
is indeed different from just projecting the unconstrained mean.

\begin{table}[tbp]
	\centering
	\caption{Algorithm comparison on the constrained Riemannian canter of mass.}
		\begin{filecontents*}{C2Dball-mean-200-benchmark.csv}
names,maxiter,times
ALM,29,0.017252758400000002
EPM,101,0.14221077480000002
PG,6,0.00027964980000000003
\end{filecontents*}
\pgfplotstabletypeset[
	col sep = comma,
	display columns/0/.style = {column name={Algorithm}, string type, column type = {r}},
	display columns/1/.style = {
		column name = Iterations,
		string type,
		column type={r},
	},
	display columns/2/.style = {
		column name = Time (msec.),
		string type,
		sci,sci zerofill, dec sep align, precision=5
	},
	every head row/.style = {
		before row = \toprule,
		after row = \midrule
	},
	every last row/.style = {after row = \bottomrule},
	]{C2Dball-mean-200-benchmark.csv}
	\label{table:SP-comparisons}
\end{table}

Statistics about the solver runs are summarized in Table~\ref{table:SP-comparisons},
where the time measurement is the mean of measuring 10 solver runs using \texttt{Chairmarks.jl}.
In these results we see that PGA needs both the least number of iterations as well as the
least run time.

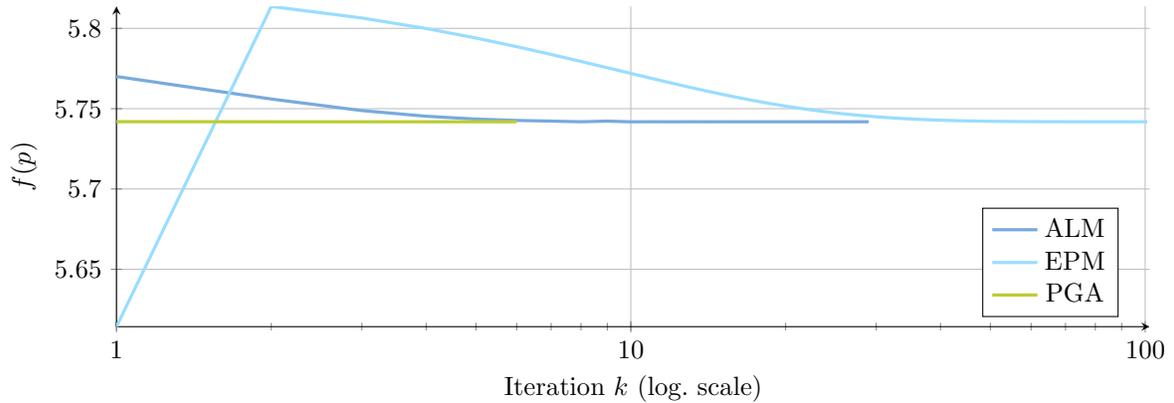
\begin{figure}
	\pgfplotstableread[col sep = comma]{C2Dball-mean-200-cost-alm.csv}\CTwoDBallCostALM
	\pgfplotstableread[col sep = comma]{C2Dball-mean-200-cost-epm.csv}\CTwoDBallCostEPM
	\pgfplotstableread[col sep = comma]{C2Dball-mean-200-cost-pg.csv}\CTwoDBallCostPG
 	\centering
 	\begin{tikzpicture}
 		\begin{axis}[
 			width = .9\textwidth,  
 			height = .25\textheight,
 			xmode = log, ymode = normal,
 			xmin = 1, xmax = 102,
 			grid = major,
 			xtick = {1e0,1e1,1e2},
 			xticklabels = {1,10,100},
 			xlabel = {Iteration $k$ (log.~scale)},
 			ylabel = {$f(p)$},
 			axis lines = left,
			legend pos=south east
 			]
 			\addplot[TolLightBlue, line width = 1.2pt] table[x = i, y = c] {\CTwoDBallCostALM};
 			\addlegendentry{ALM}
			\addplot[TolLightCyan, line width = 1.2pt] table[x = i, y = c] {\CTwoDBallCostEPM};
 			\addlegendentry{EPM}
			\addplot[TolLightPear, line width = 1.2pt] table[x = i, y = c] {\CTwoDBallCostPG};
 			\addlegendentry{PGA}
 		\end{axis}
 	\end{tikzpicture}
 	\caption{Objective versus iteration plot for different solvers on the 2D constrained mean.}
 	\label{figure:C2DBall:Cost}
\end{figure}%

When plotting the cost $f(p_k)$ during the iterations, we see in Figure~\ref{figure:C2DBall:Cost}
that the PGA already yields a very good cost in the very first step, while both
ALM with a too high cost staying inside the constrained and EPM first obtaining a lower cost
by “stepping outside” the constrained set, take much more time.

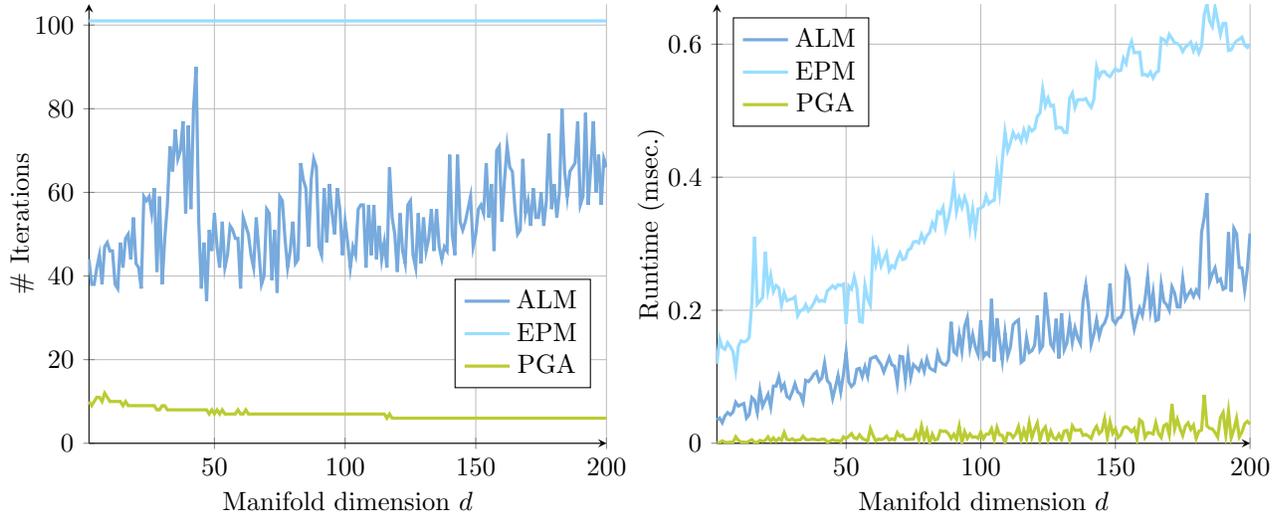
\begin{figure}
	\centering
	\pgfplotstableread[col sep = comma]{CnBallConstrMean-2-200-times.csv}\CnTimes
	\pgfplotstableread[col sep = comma]{CnBallConstrMean-2-200-iterations.csv}\CnIter
	\begin{subcaptionblock}{.49\textwidth}
		\centering
		\begin{tikzpicture}
			\begin{axis}[
				width = .825\textwidth,
				height = 0.25\textheight,
				ymin=0,
				ymax = 105,
				xmin=2, xmax=200,
				grid = major,
				axis lines = left,
				scale only axis,
				ylabel = {\# Iterations},
				xlabel = {Manifold dimension $d$},
				legend style={at={(0.97,0.25)},anchor=east},
				y label style={at={(0.05,0.5)}}
				]
				\addplot[TolLightBlue, line width = 1.2pt] table[x = d, y = alm] {\CnIter};
			\addlegendentry{ALM}
			   \addplot[TolLightCyan, line width = 1.2pt] table[x = d, y = epm] {\CnIter};
			\addlegendentry{EPM}
			   \addplot[TolLightPear, line width = 1.2pt] table[x = d, y = pga] {\CnIter};
			\addlegendentry{PGA}
			\end{axis}
		\end{tikzpicture}
	\end{subcaptionblock}%
	\begin{subcaptionblock}{.49\textwidth}
		\centering
		\begin{tikzpicture}
			\begin{axis}[
				width = .85\textwidth,
				height = 0.25\textheight,
				scale only axis,
				xmin=2, xmax=200, ymin=0,
				axis lines = left,
				grid = major,
				ylabel = {Runtime (msec.)},
				xlabel = {Manifold dimension $d$},
				y label style={at={(0.05,0.5)}},
				legend pos=north west,
				]
				\addplot[TolLightBlue, line width = 1.2pt] table[x = d, y = alm] {\CnTimes};
				\addlegendentry{ALM}
			   \addplot[TolLightCyan, line width = 1.2pt] table[x = d, y = epm] {\CnTimes};
				\addlegendentry{EPM}
			   \addplot[TolLightPear, line width = 1.2pt] table[x = d, y = pga] {\CnTimes};
				\addlegendentry{PGA}
			\end{axis}
		\end{tikzpicture}
	\end{subcaptionblock}%
	\caption{A comparison of the exact penalty method (EPM), augmented Lagrangian mathod (ALM), and the projected gradient
	for the constrained Riemannian center of mass on $\mathbb H_1^d$, $d=2,\ldots,200$:
	left: the number of iterations per experiment (right) the runtime per experimant}\label{fig:nd-mean}
\end{figure}

To continue this experiment in higher dimensions $d$\footnote{
	The complete code is available at \href{https://juliamanifolds.github.io/ManoptExamples.jl/stable/examples/Constrained-Mean-Hn/}{juliamanifolds.github.io/ManoptExamples.jl/stable/examples/Constrained-Mean-Hn/}
	}. The radius $r$ of the constrained set $C_{p_0,r}$
and the noise level $\sigma$ have to be carefully chosen, since both ALM as well as EPM otherwise
easily run out of the numeric range of \verb|Float64|.
We set for a dimension $d \in \{2,\ldots,200\}$ these values to
\begin{eqnarray}
	r_d = \frac{1}{\sqrt{n}}\quad\text{and}\quad \sigma_d = \frac{3}{2 (n-1)^{\frac{1}{4}}},
\end{eqnarray}
and choose $N=400$ random points nearly the same as before, just using
$2r_d$ times a vector of ones as the offset instead of $\frac{3}{2}(1,1)$.
Then the algorithms are set up the same way as in the previous experiment for $d=2$.
We also reuse the center point $p_0 = (0,\ldots,0,1) \in \mathbb H_1^d$ as both the
center of the constrained set as well as the start point of the algorithms.

The results of both the iterations required as well as the run time are depicted
in Fgure~\ref{fig:nd-mean}.
We see that the EPM always runs for 101 iterations, where the penality parameter
yields the algorithm to stop. These are the default parameters in \texttt{Manopt.jl}
and they maybe could be improved. ALM requires at least around 35-40 iterations, sometimes up to 77 for $d=38$.
Compared to that, the number of iterations for PGA are always below 15, for higher-dimensional
hyperbolic spaces even first 7 then only 6 iterations.

At the same time, the single iterations for PGA are relatively cheap, hence the PGA
is always fastest to finish, while ALM and even more EPM increase in there runtime.
Note that all three algorithms use in-place operations whenever possible to avoid memory allocations.
Since the number of points stays the same, this is probably due to the subproblem ALM and EPM have to solve
and corresponding memory allocations that are necessary.
Compared to that the PGA
iterations are much simpler, such that the memory allocations are less than 10\,\%
compared to one of the other two algorithms.
Together with the low number of iterations this yields that the PGA finished always faster than
the other two algorithms, in higher dimensions even in an order of magnitude.

\section{Conclusions} \label{sec:Conclusions}
In this paper, we have examined the gradient projection method for minimizing potentially non-convex functions. To further explore this topic, an intriguing
aspect that could be explored is conducting a comprehensive study on the asymptotic convergence of the method for optimization problem with  convex objective functions. It is worth noting that the inclusion of the projection in the definition of the method poses a challenging problem, and we believe that additional properties of the projection would be necessary for a successful analysis.  Given the significance of comprehending the computation of the projection in implementations of the gradient projection method, it would be worthwhile to pursue new research that offers techniques for computing the projection in additional convex sets, as well as obtaining new properties of it. This study can open up possibilities for designing other methods in hyperbolic spaces. In fact, considering that the gradient projection method in the Euclidean setting serves as a fundamental framework for the development of more sophisticated approaches, it is theoretically plausible to adapt various variants, including the Barzilai-Borwein method, spectral gradient method, and Nesterov accelerated gradient method, among others, to tackle constrained problems in hyperbolic space. In this context, we also foresee that incorporating additional properties of the projection will be crucial for successfully adapting these methods.

An essential aspect concerning the Gradient Projection Method is its numerical performance in practical applications, such as computing the constrained center of mass and assessing its efficacy compared to alternative optimization techniques. In Section~\ref{sec:Nuerics}, we present numerical experiments that validate our theoretical findings, demonstrating the practical effectiveness and efficiency of the gradient projection methods developed in this paper. It is noteworthy that variants of the Gradient Projection Method which have become popular in Euclidean settings—particularly the accelerated schemes introduced by Nesterov, such as ISTA and FISTA—are well-known for their strong performance in large-scale applications. Given this success, we anticipate similar promising results in hyperbolic spaces. Consequently, our ongoing research includes extending and adapting these accelerated gradient approaches to hyperbolic settings, which we aim to explore in forthcoming studies.


\subsection*{Funding}
O. P. Ferreira was partially supported in part by  CNPq - Brazil  Grants 304666/2021-1

\bibliographystyle{plain}
\bibliography{manuscript-arxiv}

\begin{thebibliography}{10}

\bibitem{Absil2008}
P.-A. Absil, R.~Mahony, and R.~Sepulchre.
\newblock {\em Optimization algorithms on matrix manifolds}.
\newblock Princeton University Press, Princeton, NJ, 2008.
\newblock With a foreword by Paul Van Dooren.

\bibitem{Afsari2011}
Bijan Afsari.
\newblock Riemannian {$L^p$} center of mass: existence, uniqueness, and convexity.
\newblock {\em Proc. Amer. Math. Soc.}, 139(2):655--673, 2011.

\bibitem{AfsariVidal2013}
Bijan Afsari, Roberto Tron, and Ren\'{e} Vidal.
\newblock On the convergence of gradient descent for finding the {R}iemannian center of mass.
\newblock {\em SIAM J. Control Optim.}, 51(3):2230--2260, 2013.

\bibitem{Andreanietall2007}
R.~Andreani, E.~G. Birgin, J.~M. Mart{\'{\i}}nez, and M.~L. Schuverdt.
\newblock On augmented {Lagrangian} methods with general lower-level constraints.
\newblock {\em SIAM J. Optim.}, 18(4):1286--1309, 2007.

\bibitem{AxenBaranBergmannRzecki:2023}
Seth~D. Axen, Mateusz Baran, Ronny Bergmann, and Krzysztof Rzecki.
\newblock Manifolds.jl: An extensible julia framework for data analysis on manifolds.
\newblock {\em ACM Transactions on Mathematical Software}, 49(4), 12 2023.

\bibitem{Bauschke2018}
Heinz~H. Bauschke, Minh~N. Bui, and Xianfu Wang.
\newblock Projecting onto the intersection of a cone and a sphere.
\newblock {\em SIAM J. Optim.}, 28(3):2158--2188, 2018.

\bibitem{Bacak2014}
Miroslav Ba\v{c}\'{a}k.
\newblock Computing medians and means in {H}adamard spaces.
\newblock {\em SIAM J. Optim.}, 24(3):1542--1566, 2014.

\bibitem{BenedettiPetronio1992}
Riccardo Benedetti and Carlo Petronio.
\newblock {\em Lectures on hyperbolic geometry}.
\newblock Universitext. Springer-Verlag, Berlin, 1992.

\bibitem{Bergmann:2022}
Ronny Bergmann.
\newblock Manopt.jl: Optimization on manifolds in {J}ulia.
\newblock {\em Journal of Open Source Software}, 7(70):3866, 2022.

\bibitem{Boumal2020}
Nicolas Boumal.
\newblock {\em An introduction to optimization on smooth manifolds}.
\newblock Cambridge University Press, Cambridge, 2023.

\bibitem{Cannon1997}
James~W. Cannon, William~J. Floyd, Richard Kenyon, and Walter~R. Parry.
\newblock Hyperbolic geometry.
\newblock In {\em Flavors of geometry}, volume~31 of {\em Math. Sci. Res. Inst. Publ.}, pages 59--115. Cambridge Univ. Press, Cambridge, 1997.

\bibitem{chami2021horopca}
Ines Chami, Albert Gu, Dat~P Nguyen, and Christopher R{\'e}.
\newblock Horopca: Hyperbolic dimensionality reduction via horospherical projections.
\newblock In {\em International Conference on Machine Learning}, pages 1419--1429. PMLR, 2021.

\bibitem{Chen2010}
Donghui Chen and Robert~J. Plemmons.
\newblock Nonnegativity constraints in numerical analysis.
\newblock In {\em The birth of numerical analysis}, pages 109--139. World Sci. Publ., Hackensack, NJ, 2010.

\bibitem{criscitiello2022negative}
Christopher Criscitiello and Nicolas Boumal.
\newblock Negative curvature obstructs acceleration for strongly geodesically convex optimization, even with exact first-order oracles.
\newblock In {\em Conference on Learning Theory}, pages 496--542. PMLR, 2022.

\bibitem{Faber1983}
Richard~L. Faber.
\newblock {\em Differential geometry and relativity theory. {An} introduction}, volume~76 of {\em Pure Appl. Math., Marcel Dekker}.
\newblock Marcel Dekker, Inc., New York, NY, 1983.

\bibitem{FacchineiPangI}
Francisco Facchinei and Jong-Shi Pang.
\newblock {\em Finite-dimensional variational inequalities and complementarity problems. {V}ol. {I}}.
\newblock Springer Series in Operations Research. Springer-Verlag, New York, 2003.

\bibitem{fan2022nested}
Xiran Fan, Chun-Hao Yang, and Baba~C Vemuri.
\newblock Nested hyperbolic spaces for dimensionality reduction and hyperbolic nn design.
\newblock In {\em Proceedings of the IEEE/CVF Conference on Computer Vision and Pattern Recognition}, pages 356--365, 2022.

\bibitem{Feng2022}
Shuailing Feng, Wen Huang, Lele Song, Shihui Ying, and Tieyong Zeng.
\newblock Proximal gradient method for nonconvex and nonsmooth optimization on {H}adamard manifolds.
\newblock {\em Optim. Lett.}, 16(8):2277--2297, 2022.

\bibitem{FerreiraIusemSandor2013}
O.~P. Ferreira, A.~N. Iusem, and S.~Z. N\'{e}meth.
\newblock Projections onto convex sets on the sphere.
\newblock {\em J. Global Optim.}, 57(3):663--676, 2013.

\bibitem{FerreiraOliveira2002}
O.~P. Ferreira and P.~R. Oliveira.
\newblock Proximal point algorithm on {R}iemannian manifolds.
\newblock {\em Optimization}, 51(2):257--270, 2002.

\bibitem{FerreiraNemethShu2022}
Orizon~Pereira Ferreira, S{\'a}ndor~Zolt{\'a}n N{\'e}meth, and Jinzhen Zhu.
\newblock Convexity of sets and functions on the hyperbolic space.
\newblock {\em J Optim Theory Appl}, 2022.

\bibitem{fletcher2004principal}
P~Thomas Fletcher, Conglin Lu, Stephen~M Pizer, and Sarang Joshi.
\newblock Principal geodesic analysis for the study of nonlinear statistics of shape.
\newblock {\em IEEE transactions on medical imaging}, 23(8):995--1005, 2004.

\bibitem{HiriartLemarechal2001}
Jean-Baptiste Hiriart-Urruty and Claude Lemar\'{e}chal.
\newblock {\em Fundamentals of convex analysis}.
\newblock Grundlehren Text Editions. Springer-Verlag, Berlin, 2001.
\newblock Abridged version of {{\i}t Convex analysis and minimization algorithms. I} [Springer, Berlin, 1993; MR1261420 (95m:90001)] and {{\i}t II} [ibid.; MR1295240 (95m:90002)].

\bibitem{HuangWei2022}
Wen Huang and Ke~Wei.
\newblock Riemannian proximal gradient methods.
\newblock {\em Math. Program.}, 194(1-2):371--413, 2022.

\bibitem{jawanpuria2019low}
Pratik Jawanpuria, Mayank Meghwanshi, and Bamdev Mishra.
\newblock Low-rank approximations of hyperbolic embeddings.
\newblock In {\em 2019 IEEE 58th Conference on Decision and Control (CDC)}, pages 7159--7164. IEEE, 2019.

\bibitem{Karcher:1977}
H.~Karcher.
\newblock Riemannian center of mass and mollifier smoothing.
\newblock {\em Communications on Pure and Applied Mathematics}, 30(5):509--541, 1977.

\bibitem{keller2021hyperbolic}
Martin Keller-Ressel and Stephanie Nargang.
\newblock The hyperbolic geometry of financial networks.
\newblock {\em Scientific reports}, 11(1):1--12, 2021.

\bibitem{Kriouko2010}
Dmitri Krioukov, Fragkiskos Papadopoulos, Maksim Kitsak, Amin Vahdat, and Mari\'{a}n Bogu\~{n}\'{a}.
\newblock Hyperbolic geometry of complex networks.
\newblock {\em Phys. Rev. E (3)}, 82(3):036106, 18, 2010.

\bibitem{LaiYoshise2024}
Zhijian Lai and Akiko Yoshise.
\newblock Riemannian interior point methods for constrained optimization on manifolds.
\newblock {\em J. Optim. Theory Appl.}, 201(1):433--469, 2024.

\bibitem{Lewis2022}
Adrian~S. Lewis, Genaro Lopez-Acedo, and Adriana Nicolae.
\newblock Local linear convergence of alternating projections in metric spaces with bounded curvature.
\newblock {\em SIAM J. Optim.}, 32(2):1094--1119, 2022.

\bibitem{LiuBoumal2020}
Changshuo Liu and Nicolas Boumal.
\newblock Simple algorithms for optimization on {Riemannian} manifolds with constraints.
\newblock {\em Appl. Math. Optim.}, 82(3):949--981, 2020.

\bibitem{pmlr-v119-lou20a}
Aaron Lou, Isay Katsman, Qingxuan Jiang, Serge Belongie, Ser-Nam Lim, and Christopher De~Sa.
\newblock Differentiating through the fr{é}chet mean.
\newblock In Hal~Daumé III and Aarti Singh, editors, {\em Proceedings of the 37th International Conference on Machine Learning}, volume 119 of {\em Proceedings of Machine Learning Research}, pages 6393--6403. PMLR, 13--18 Jul 2020.

\bibitem{Martinez023}
David Martinez-Rubio, Christophe Roux, Christopher Criscitiello, and Sebastian Pokutta.
\newblock Accelerated methods for riemannian min-max optimization ensuring bounded geometric penalties.
\newblock {\em arXiv preprint arXiv:2305.16186}, 2023.

\bibitem{Moreira2024}
Gabriel Moreira, Manuel Marques, Jo{\~a}o~Paulo Costeira, and Alexander Hauptmann.
\newblock Hyperbolic vs euclidean embeddings in few-shot learning: Two sides of the same coin.
\newblock In {\em Proceedings of the IEEE/CVF Winter Conference on Applications of Computer Vision}, pages 2082--2090, 2024.

\bibitem{Moshir2021}
Mahdi Moshiri, Farshad Safaei, and Zeynab Samei.
\newblock A novel recovery strategy based on link prediction and hyperbolic geometry of complex networks.
\newblock {\em J. Complex Netw.}, 9(4):Paper No. cnab007, 2021.

\bibitem{muscoloni2017machine}
Alessandro Muscoloni, Josephine~Maria Thomas, Sara Ciucci, Ginestra Bianconi, and Carlo~Vittorio Cannistraci.
\newblock Machine learning meets complex networks via coalescent embedding in the hyperbolic space.
\newblock {\em Nature communications}, 8(1):1--19, 2017.

\bibitem{nickel2018learning}
Maximillian Nickel and Douwe Kiela.
\newblock Learning continuous hierarchies in the lorentz model of hyperbolic geometry.
\newblock In {\em International Conference on Machine Learning}, pages 3779--3788. PMLR, 2018.

\bibitem{Obara2022}
Mitsuaki Obara, Takayuki Okuno, and Akiko Takeda.
\newblock Sequential quadratic optimization for nonlinear optimization problems on {Riemannian} manifolds.
\newblock {\em SIAM J. Optim.}, 32(2):822--853, 2022.

\bibitem{peng2021hyperbolic}
Wei Peng, Tuomas Varanka, Abdelrahman Mostafa, Henglin Shi, and Guoying Zhao.
\newblock Hyperbolic deep neural networks: A survey.
\newblock {\em IEEE Transactions on Pattern Analysis and Machine Intelligence}, 44(12):10023--10044, 2021.

\bibitem{Ratcliffe2019}
John~G. Ratcliffe.
\newblock {\em Foundations of hyperbolic manifolds}, volume 149 of {\em Graduate Texts in Mathematics}.
\newblock Springer, Cham, third edition, [2019] \copyright 2019.

\bibitem{Sakai1996}
Takashi Sakai.
\newblock {\em Riemannian geometry}, volume 149 of {\em Translations of Mathematical Monographs}.
\newblock American Mathematical Society, Providence, RI, 1996.
\newblock Translated from the 1992 Japanese original by the author.

\bibitem{sharpee2019argument}
Tatyana~O Sharpee.
\newblock An argument for hyperbolic geometry in neural circuits.
\newblock {\em Current opinion in neurobiology}, 58:101--104, 2019.

\bibitem{Tabaghi2021}
Puoya Tabaghi and Ivan Dokmani.
\newblock On procrustes analysis in hyperbolic space.
\newblock {\em IEEE Signal Processing Letters}, 28:1120--1124, 2021.

\bibitem{tabaghi2020hyperbolic}
Puoya Tabaghi and Ivan Dokmani{\'c}.
\newblock Hyperbolic distance matrices.
\newblock In {\em Proceedings of the 26th ACM SIGKDD International Conference on Knowledge Discovery \& Data Mining}, pages 1728--1738, 2020.

\bibitem{Walter1974}
Rolf Walter.
\newblock On the metric projection onto convex sets in riemannian spaces.
\newblock {\em Arch. Math. (Basel)}, 25:91--98, 1974.

\bibitem{WeberSra2023}
Melanie Weber and Suvrit Sra.
\newblock Riemannian optimization via {Frank}-{Wolfe} methods.
\newblock {\em Math. Program.}, 199(1-2 (A)):525--556, 2023.

\bibitem{wilson2018gradient}
Benjamin Wilson and Matthias Leimeister.
\newblock Gradient descent in hyperbolic space, 2018.

\bibitem{Wilson2014}
Richard~C. Wilson, Edwin~R. Hancock, El?bieta Pekalska, and Robert~P.W. Duin.
\newblock Spherical and hyperbolic embeddings of data.
\newblock {\em IEEE Transactions on Pattern Analysis and Machine Intelligence}, 36(11):2255--2269, 2014.

\bibitem{zhang2016}
Hongyi Zhang and Suvrit Sra.
\newblock First-order methods for geodesically convex optimization.
\newblock In {\em Conference on Learning Theory}, pages 1617--1638. PMLR, 2016.

\end{thebibliography}


\end{document}